\documentclass{article}
\usepackage[T2A]{fontenc}
\usepackage[cp1251]{inputenc}
\usepackage[tbtags]{amsmath}
\usepackage{amsfonts,amssymb,mathrsfs,amscd}
\usepackage{graphicx}

\def\thtext#1{
  \catcode`@=11
  \gdef\@thmcountersep{. #1}
  \catcode`@=12
}

\def\threst{
  \catcode`@=11
  \gdef\@thmcountersep{.}
  \catcode`@=12
}

\newtheorem{thm}{Theorem}[section]

\newtheorem{ass}{Assertion}[section]
\newtheorem{lem}{Lemma}[section]
\newtheorem{cor}{Corollary}[section]
\newtheorem{conj}{Conjecture}[section]
\newtheorem{prb}{Problem}[section]

{\endtrivlist}
\newenvironment{examp}{\trivlist \item[\hskip \labelsep{\bf Example.}]}%
{\endtrivlist}

\newenvironment{rk}{\trivlist \item[\hskip \labelsep{\bf Remark.}]}%
{\endtrivlist}
{\endtrivlist}
\newenvironment{proof}{\trivlist \item[\hskip \labelsep{\bf
Proof.}]}{\endtrivlist}
{\endtrivlist}

\newdimen\fgm
 \def\Pic#1#2#3#4{
   \begin{figure}
\begin{center}\includegraphics[height=#4\fgm]{#1}\end{center}
     \caption{\label{#3} #2}
  \end{figure}
 }

\catcode`@=11

\def\.{.\spacefactor\@m}

\catcode`@=12

\def\<{\langle}
\def\>{\rangle}
\def\r{\rho}
\def\v{\varphi}
\def\:{\colon}
\def\om{\omega}
\def\R{{\Bbb R}}

\def\g{\gamma}
\def\dl{\delta}
\def\ss{\subset}
\def\d{\partial}
\def\sm{\setminus}
\def\bG{{\bar G}}

\def\bom{{\bar\om}}

\def\mst{\operatorname{mst}}
\def\cM{{\cal M}}
\def\e{\varepsilon}
\def\bp{{\bar p}}

\def\s{\sigma}
\def\rom#1{{\em #1}}
\def\({{\rom(}}
\def\){{\rom)}}

\def\ssr{\operatorname{ssr}}

\def\smt{\operatorname{smt}}
\def\mf{\operatorname{mf}}
\def\c{\circ}

\def\a{\alpha}
\def\b{\beta}

\def\G{{\Gamma}}
\def\cX{{\cal X}}
\def\mpn{\operatorname{mpn}}

\def\0{\emptyset}
\def\cW{{\cal W}}
\def\cG{{\cal G}}
\def\rom#1{{\em #1}}
\def\({{\rom(}}
\def\){{\rom)}}

\def\l{{\lambda}}
\def\cO{{\cal O}}
\def\sgr{\operatorname{sgr}}

\def\cI{{\cal I}}
\def\cE{{\cal E}}
\def\cH{{\cal H}}
\def\cN{{\cal N}}
\def\cL{{\cal L}}
\def\cY{{\cal Y}}

\def\Om{\Omega}
\def\mpf{\operatorname{mpf}}
\def\Q{{\Bbb Q}}
\def\bE{{\bar E}}

\def\degen{\operatorname{degen}}
\def\cP{{\cal P}}
\def\Y{{\Upsilon}}
\def\neind{\operatorname{neind}}
\def\bg{{\bar\g}}
\def\ti{\widetilde}
\def\sr{\operatorname{sr}}
\begin{document}

\title{One-dimensional Gromov minimal filling}
\author{A.\,O.~Ivanov and A.\,A.~Tuzhilin}

\maketitle

 \begin{abstract}
The present paper opens a new branch in the theory of variational problems with branching extremals, the investigation of one-dimensional minimal fillings of finite pseudo-metric spaces. On the one hand, this problem is a one-dimensional version of a generalization of Gromov's minimal fillings problem to the case of stratified manifolds (the filling in our case is a weighted graph). On the other hand, this problem is interesting in itself and also can be considered as a generalization of another classical problem, namely, the Steiner problem on the construction of a shortest network joining a given set of terminals. Besides the statement of the problem, we discuss several properties of the minimal fillings, describe minimal fillings of additive spaces, and state several conjectures. We also include some announcements concerning the very recent results obtained in our group, including a formula calculating the weight of the minimal filling for an arbitrary finite pseudo-metric space and the concept of pseudo-additive space which generalizes the classical concept of additive space. We hope that the theory of one-dimensional minimal fillings refreshes the interest in the Steiner problem and gives an opportunity to solve several long standing problems, such as the calculation of the Steiner ratio, in particular the verification of the Gilbert--Pollack conjecture on the Steiner ratio of the Euclidean plane.

Bibliography: 33 items.
 \end{abstract}


\markright{One-dimensional Gromov minimal fillings problem}

 \section*{Introduction}
 \markright {Introduction.}
The problem considered in this paper appears as a result of a synthesis of two classical problems: the Steiner problem on the shortest networks, and Gromov's problem on minimal fillings. Recall briefly the history of these questions.

The Steiner problem is a problem on optimal connection of a finite set of points of a metric space. Apparently, the first problems of that type appeared in works of P.~Fermat, who stated the question on finding the location of a point in the plane, such that the sum of the distances from it to the vertices of a given triangle is minimal. It took a few centuries to obtain a complete answer (Torricelli, Simpson, etc., see details in~\cite{ITRFFIBook}), which demonstrates that joining the points in the plane it might be profitable to add an additional point-fork. The importance of such additional points was quite clear to C.~F.~Gauss discussing in his letters to Schuhmacher the construction of the shortest road network joining famous German cities Hamburg, Bremen, Hannower and Braunschweig. In 1934 Jarnik and K\"ossler~\cite{JaKo} stated the general problem which is now known as Classical Steiner Problem. In fact, their problem generalizes the problems of Fermat and Gauss on the shortest connection to the case of an arbitrary finite set of points in the plane. As concerns Steiner, he works with another generalization of the Fermat problem: to find  a point in the space such that the sum of the distances from it to the given ones is minimal. Notice that the priority misunderstanding appears due to popular book of Courant and Robbins ``What is Mathematics?''~\cite{CourRob}, where Fermat problem is referred as Steiner problem, and the Jarnik and K\"ossler problem is called just a generalization of the Steiner problem.

There are a lot of publications devoted to the classical Steiner problem and its numerous generalizations, for example, to the cases of normed spaces, Riemannian manifolds, and more general metric spaces, say, spaces of words. An interested reader could find a review in~\cite{ITRFFIBook}. Another surge of interest in Steiner problem is connected to Gilbert--Pollack Conjecture~\cite{GilPol} on the Steiner ratio of Euclidean plane. Numerous attempts to prove it have failed. The best known attempt belongs to Du and Hwang~\cite{DuHwang}, but it turns out that their reasoning contains serious gaps, see~\cite{ITLup},~\cite{ITRFFIBook},~\cite{deWet1},~\cite{deWet2},~\cite{IKMS}. Ideas of Du and Hwang were advertised in the stage of announcing publications~\cite{NYT} that has led to popularization of the wrong construction.  Several papers appeared (for example, \cite{Rubin},~\cite{Lambda}) where the ideas of Du and Hwang were adopted to the case under consideration, and, as a result, some unfounded conclusions got the status of theorems. Notice that the validity of Gilbert--Pollack conjecture itself seems undoubted, therefore attempts to prove it appear again and again. In particular, numerous authors, including the authors of the present paper, tried to improve the construction of Du and Hwang, but did not succeed. It might make sense to search for a completely different approach to the problem. The minimal fillings discussed in the present paper could be a base for such an approach.

The concept of a minimal filling appeared in papers of Gromov~\cite{Gromov} in the following form. Let $M$ be a smooth closed manifold endowed with a distance function  $\r$. Consider all possible films $W$ spanning $M$, i.e. compact manifolds with the boundary $M$. Consider on  $W$ a distance function $d$ non-decreasing the distances between the points from  $M$. Such a metric space $\cW=(W,d)$ is called by a {\em filling\/} of the metric space $\cM=(M,\r)$. The Gromov Problem consists in describing the infimum of the volumes of the fillings and in searching for the spaces $\cW$ where this infimum is achieved and which are called {\em minimal fillings}.

An interest in minimal fillings is inspired, first of all, by the fact that many classical geometrical inequalities can be stated in terms of the fillings; see details and exact references in the monograph~\cite{BurBurIva} and also in the dissertation of S.~V.~Ivanov~\cite{SIvanovDis}. For example, the Bezikovich inequality, which says that the volume of a Riemannian cube is greater than or equal to the product of the distances between its opposite faces, follows from the fact that the standard cube in the Euclidean space is the minimal filling of its boundary. Another example is the Pu inequality giving a lower estimate for the length of the shortest non-contractible loop on the projective plane. This inequality follows from the fact that the standard half-sphere is a minimal filling of its boundary circle endowed with the intrinsic metric. Notice also that minimal fillings possess numerous applications in dynamic systems theory, asymptotic geometry, mathematical physics, etc.

In the scope of Steiner problem, it is natural to consider $M$ as a finite metric space. Then the possible fillings are metric spaces having the structure of one-dimensional stratified manifolds which can be considered as edge-weighted graphs with non-negative weight function. This leads to the following particular case of generalized Gromov problem.

Let $M$ be an arbitrary finite set, and $G=(V,E)$ be a connected graph. We say, that $G$ {\em joins $M$}, if $M\ss V$. Now, let $\cM=(M,\r)$ be a finite pseudo-metric space (distances between the distinct points can be equal to zero), $G=(V,E)$ be a connected graph joining $M$, and $\om\:E\to\R_+$ is a mapping into non-negative numbers, which is usually referred as a {\em weight function} and which generates a {\em weighted graph\/} $\cG=(G,\om)$. The function $\om$ generates on $V$ the pseudo-metric $d_\om$, namely, the distance between the vertices of the graph  $\cG$ is defined as the least possible weight of the path in $\cG$ joining these vertices. If for any two points $p$ and $q$ from $M$ the inequality $\r(p,q)\le d_\om(p,q)$ is valid, then the weighted graph $\cG$ is called by a {\em filing\/} of the space $\cM$, and the graph $G$ is referred as the {\em type\/} of this filing. The value $\mf(\cM)=\inf\om(\cG)$, where infimum is taken over all the filings $\cG$ of the space $\cM$ is called by the {\em weight of minimal filling}, and a filling $\cG$ such that  $\om(\cG)=\mf(\cM)$ is called by a {\em minimal filling}. The main problem is to learn how to calculate $\mf(\cM)$ and how to describe minimal fillings.

Notice that stratified manifolds have appeared naturally in geometric problems; see for example~\cite{FomReif}, \cite{FomPlat}, \cite{Gusev}, and in such applications as quantum physics~\cite{ITDeformStr}, \cite{Shafar}.

Let us describe briefly the content and the main results of the paper. In Preliminaries we introduce the main objects of our investigation and discuss some relations between them. In general these objects may have very complicated topologies, however, one can essentially simplify the structure without loss of the sense of the problem. Section~\ref{sec:reduction}, Reduction to trees, is devoted to such simplifications. It turns out that in the case of shortest networks and minimal fillings one can consider only trees without interior vertices of degree $2$; moreover, if we alow the degenerate edges (the edges of zero weight or length) then it is sufficient to use only binary trees, i.e. the trees with vertices of degree $1$ and $3$ only, see Theorem~\ref{th:trees_fill}. Section~\ref{sec:realization}, Minimal realization, shows that each minimal filling can be represented as a shortest tree with some special boundary in some special metric space. Namely, see Corollaries~\ref{cor:Kurat_mpf} and~\ref{cor:Kur_mf}, if $\cM=(M,\r)$ is a finite pseudo-metric space consisting of $n$ points, then the Kuratowski isometrical embedding $\v_\cM\:M\to\ell_\infty^n$ can be extended to an isometrical embedding of a minimal (parametric) filling $\cG=(G,\om)$, such that the image $\v_\cM(G)$ is a shortest tree (a minimal parametric tree) joining the boundary $\v_\cM(M)\ss\ell_\infty^n$. In Section~\ref{sec:reformul} we discuss the relations between the concepts of minimal filling, shortest tree, and minimal spanning tree. It is shown that the weight of minimal filling can be found as infimum of the lengths of minimal spanning trees for so-called extensions of the initial finite metric space, see Corollaries~\ref{cor:secondDefOfMF} and~\ref{cor:thirdDefOfMF}. In Section~\ref{sec:exist} we show that each minimal parametric filling is a solution to a linear programming problem. Using this idea, we prove that for any finite metric space the set of all its minimal parametric fillings of a fixed topology $G$ is a non-empty convex compact in the configuration space of all weighted trees $(G,\om)$, see Assertion~\ref{ass:opt_weight}.

Further, in Section~\ref{sec:exact} the concept of an exact path in a filling is introduced. By definition of filling, for any two points $p$ and $q$ of a finite metric space $\cM=(M,\r)$ and any its filling $\cG=(G,\om)$, the inequality $\om(\g_{pq})\ge\r(p,q)$ is valid, where $\g_{pq}$ is the path in $G$ joining $p$ and $q$. Such a boundary path $\g_{pq}$ is said to be {\em exact}, if $\om(\g_{pq})=\r(p,q)$, and a path (an edge, a vertex) is said to be exact if it is contained in an exact boundary path. Section~\ref{sec:exact} describes the structure of the set of exact paths in minimal (parametric) fillings. In the subsequent sections this technique gets several applications. First such application appears in Section~\ref{sec:tours}. For any cyclic order $\pi\:M\to M$ on $M$ we define the corresponding {\em half-perimeter\/} $p(\cM,\pi)$ of the pseudo-metric space $\cM=(M,\r)$ as follows:
 $$
p(\cM,\pi)=\frac12\sum_{p\in M}\r\bigl(p,\pi(p)\bigr).
 $$
Let $G$ be a tree joining $M$. A cyclic order is called a {\em tour around $G$}, if it corresponds to an Euler tour of the doubling of $G$. Then, see Corollary~\ref{cor:fill_ge_halfp}, for any filling $\cG=(G,\om)$ of a finite pseudo-metric space $\cM=(M,\r)$, we prove that
 $$
\om(G)\ge \max_{\pi\in\cO(G)} p(\cM,\pi),
 $$
where $\cO(G)$ stands for the set of all the tours around $G$. Hence, if the equality
 $$
\om(G)=\max_{\pi\in\cO(G)} p(\cM,\pi)
 $$
holds, then $\cG$ is a minimal parametric filling, and
 $$
\mpf(\cM,G)=\max_{\pi\in\cO(G)} p(\cM,\pi). \eqno{(\dagger)}
 $$
We prove, see Assertion~\ref{ass:exact_order}, that the tour $\pi$ where this maximum is achieved is exact, i.e., the boundary paths joining the points of $\cM$ consecutive with respect to $\pi$, are exact. Vice versa, if a tour $\pi$ around $G$ is exact, then $\cG=(G,\om)$ is a minimal parametric filling and the equality holds. But it is not difficult to construct an example of a minimal parametric filling with degenerate non-exact edges and without exact tours. Moreover, there exist non-degenerate minimal parametric fillings which do not possess exact tours, see Example of Z.~Ovsjannikov in Section~\ref{sec:tours}. Hence, Formula~$(\dagger)$ is not valid in general case. However, the minimal parametric fillings in the above examples are not minimal fillings, and the authors conjectured, see Conjecture~\ref{conj:min-fill-exact-tour}, that each minimal filling must possess an exact tour and, moreover, see Conjecture~\ref{conj:min-fill-formula}, the weight of minimal filling can be calculated by the next formula:
 $$
\mf(\cM)=\min_G\max_{\pi\in\cO(G)} p(\cM,\pi). \eqno{(\ddagger)}
 $$

Recent results of our group, see~\cite{IOST}, \cite{Eremin}, \cite{Ovs}, and Sections~\ref{sec:formula} and~\ref{sec:additive}, show that sometimes it is usefull to allow the weight function of a filling to take negative values. The corresponding objects defined in~\cite{IOST} are referred as {\em generalized fillings}. It turns out, that among minimal generalized fillings of an arbitrary pseudo-metric space there exists a minimal filling, see Theorem~\ref{th:IOST}. Therefore one can calculate the weight of minimal generalized filling instead of ordinary minimal filling. Investigating  the generalized fillings, A.~Eremin has found out that
Conjectures~\ref{conj:min-fill-exact-tour} and~\ref{conj:min-fill-formula} are not valid, but Formula~$(\ddagger)$ can be turned to a true formula, see Theorem~\ref{th:eremin}, by means of the generalized fillings and the concept of a multi-tour introduced by Eremin in~\cite{Eremin}. A {\em multi-tour of multiplicity\/} $k$ around a tree $G$ joining $M$ can be defined as a set of boundary paths forming an Euler tour in the $2k$-plication of $G$ (that is, the multigraph obtained from $G$ by taking $2k$ copies of each its edge). The half-perimeter of a multi-tour is defined as the sum of the corresponding $kn$ distances divided by $2k$. Eremin proved that even Formula~$(\dagger)$ becomes true providing its left-hand part means the weight of minimal parametric generalized filling, and the maximum in the right-hand part is taken over all multi-tours around $G$.

Section~\ref{sec:additive} is devoted to the additive spaces that appear in many applications. Recall that a pseudo-metric space $\cM=(M,\r)$ is additive, if there exists a weighted tree $\cG=(G,\om)$ joining $M$, such that $d_\om=\r$.  In this case the tree $\cG$ is said to be {\em generating}. A well-known description of additive spaces is given in terms of so-called {\em $4$-points rule\/}: the space is additive, if and only if for any
four points $p_i$, $p_j$, $p_k$, $p_l$, the values $\r(p_i,p_j)+\r(p_k,p_l)$, $\r(p_i,p_k)+\r(p_j,p_l)$, $\r(p_i,p_l)+\r(p_j,p_k)$ are the lengths of sides of an isosceles triangle whose base does not exceed its other sides. Effective algorithms restoring a generating tree are well-known also. For additive spaces we obtain a complete description of minimal filling, see Theorem~\ref{th:additive=minimum}, namely, we show that the set of minimal fillings of an additive space coincides with the set of its generating trees. Moreover, since a non-degenerate generating tree of an additive space is unique, then a non-degenerate minimal filling of an additive space is unique as well, see Corollary~\ref{cor:rigid}.

Further, the {\em half-perimeter of a pseudo-metric space $\cM=(M,\r)$} is defined as the value $p(\cM)=\min_\s p(\cM,\s)$, where minimum is taken over all cyclic orders on $M$. For an additive space, the weight of its minimal filling is equal to the half-perimeter of the space and is also equal to the half-perimeter of the minimal filling with respect to each its tour, see Corollary~\ref{cor:weight-mf-additive}. Moreover, O.~Rubleva shows, see~\cite{Rubleva} and Assertion~\ref{ass:Rubleva}, that a pseudo-metric space is additive, if and only if the weight of its minimal filling is equal to the half-perimeter of the space.

Since all the half-perimeters corresponding to a minimal filling of an additive space are the same, the authors stated the problem to describe all the spaces $\cM=(M,\r)$ such that all the half-perimeters corresponding to the tours around some tree joining $M$ are the same. The complete answer is obtained by
Z.~Ovsjannikov~\cite{Ovs}. He shows that the distance function of such a space $\cM=(M,\r)$ is generated by a weighted tree with possibly negative weights of some edges (such tree is also said to be {\em generating\/}). Z.~Ovsjannikov suggests to call such spaces {\em pseudo-additive\/} and shows that each such space possesses the above property. Similarly to additive spaces, the pseudo-additive spaces can be described in terms of so-called {\em weak $4$-points rule}: the space is pseudo-additive, if and only if for any four points $p_i$, $p_j$, $p_k$, $p_l$, the values $\r(p_i,p_j)+\r(p_k,p_l)$, $\r(p_i,p_k)+\r(p_j,p_l)$,  $\r(p_i,p_l)+\r(p_j,p_k)$ are the lengths of sides of an isosceles triangle, see~\cite{Ovs}. There exist effective algorithms restoring a generating tree, similar to classical ones.

In Section~\ref{sec:rays} we consider families of pseudo-metric spaces referred as rays. The {\em multiplicative ray\/} generated by a pseudo-metric space $\cM=(M,\r)$ consists of all the spaces of the form $(M,\l\r)$, $\l\ge0$. The {\em additive ray\/} generated by $\cM=(M,\r)$ consists of all the spaces of the form $(M,\r+a)$, $a\ge a_\cM$, where $a_\cM$ is a non-positive constant depending on $\cM$. Theorem~\ref{th:rays} describes minimal fillings of the elements of the rays in terms of minimal fillings of the initial space. In particular, it is shown that all the spaces from the family $(M,\l\r+a)$, $\l>0$, $a>\l a_\cM$, have the same set of types of minimal fillings.

In Section~\ref{sec:incommens} we consider finite pseudo-metric spaces, all whose non-zero distances are linearly independent over the field $\Q$ of rational numbers. We call such spaces {\em incommensurable}. It is shown that any minimal filling of such space consisting of more than $3$ points can not be ``very degenerate'', namely, it must contain a non-degenerate interior edge. Moreover, boundary edges of any minimal filling of an incommensurable metric space are also non-degenerate(see Section~\ref{sec:exact}).

In Section~\ref{sec:examp} we give several examples of minimal filling finding. In particular, it contains complete answers for three-points and four-points pseudo-metric spaces, for any regular simplex (a metric space with constant non-zero distance function), etc.

Section~\ref{sec:ratios} is inspired by Steiner ratio investigations, see above. Recall that the {\em Steiner ratio of a finite subset\/} $M$ of a metric space $\cX=(X,\r)$ is defined as
 $$
\sr(M)=\frac{\smt(M,\r)}{\mst(M,\r)},
 $$
where $\mst(M,\r)$ and $\smt(M,\r)$ stand for the lengths of minimal spanning tree and Steiner minimal tree of $M\ss X$, respectively. Then the {\em Steiner ratio $\sr\cX$ of the metric space $\cX$} and the {\em degree $n$ Steiner ratio $\sr_n\cX$} of the space $\cX$ can be defined as follows
 $$
\sr_n\cX=\inf\bigl\{\sr(M):M\ss X,\ |M|\le n\bigr\}, \qquad \sr(\cX)=\inf_{n\ge2}\sr_n\cX.
 $$
We suggest to consider the values
 $$
\sgr(M)=\frac{\mf(M,\r)}{\mst(M,\r)},\qquad \ssr(M)=\frac{\mf(M,\r)}{\smt(M,\r)},
 $$
which we call the {\em Steiner--Gromov ratio\/} and the {\em Steiner subratio\/} of $M\ss X$, respectively. Similarly to Steiner ratio case, one can define the values $\sgr_n\cX$, $\sgr\cX$, $\ssr_n\cX$, and $\ssr\cX$. We hope that the calculation of these values can be simpler than the one of classical Steiner ratio, an that there are some connections between the three values $\sgr\cX$, $\ssr\cX$, and $\sr\cX$. Section~\ref{sec:ratios} contains several preliminary results concerning the ratios. It is shown that $1/2\le\sgr(\cX)\le1$ for any pseudo-metric space and there exists a metric space whose Steiner--Gromov ratio equals $1/2$. Also it is shown, that if a pseudo-metric space $\cX$ contains an equilateral triangle, then   $\sgr_3\cX=3/4$. For the Steiner subratio, it is shown that $\ssr_3(\R^n)=\sqrt3/2$ and $\ssr_4(\R^2)=\sqrt3/2$, see E.~Filonenko~\cite{Filonenko}. We conjecture, see Conjecture~\ref{conj:subrat}, that the Steiner subratio of the Euclidean plane is equal to $\sqrt3/2$ and is achieved at the vertex set of an equilateral triangle.

The authors use the opportunity to express their heartfelt gratitude to A.~T.~Fomenko for his longstanding attention to our work, to N.~P.~Strelkova for attentive reading of the manuscript and several useful remarks, to F.~Morgan for many remarks improved essentially the quality of our English translation, and to all the participants of the seminar ``Minimal Networks'' which the authors led at mechanical and mathematical faculty of Moscow State University for their interest and numerous fruitful discussions.

 \section {Preliminaries} \label{sec:prelim}
 \markright {\thesection.~Preliminaries.}
In the present paragraph we discuss a few more optimization problems closely related with description of minimal fillings. In particular, we give formal definition of the shortest networks.

Let $\cX=(X,d)$ be a pseudo-metric space and $G=(V,E)$ an arbitrary connected graph. Any mapping $\G\:V\to X$ is called a {\em network in
$\cX$ parameterized by the graph $G=(V,E)$}, or a {\em network of the type $G$}. The {\em vertices} and {\em edges\/} of the network $\G$ are the restrictions of the mapping $\G$, respectively on the vertices and edges of the graph $G$. The {\em length of the edge\/}
$\G\:vw\to X$ is the value $d\bigl(\G(v),\G(w)\bigr)$, and the {\em length $d(\G)$ of the network $\G$} is the sum of lengths of all its edges. In what follows we shall consider various boundary value problems for graphs. To do that, we fix some subsets $\d G$ of the vertex sets $V$ of our graphs $G=(V,E)$, and we call such $\d G$ the {\em boundaries}. We always suppose that in each graph under consideration there was chosen a boundary, possibly, an empty one. The {\em boundary $\d\G$ of a network $\G$} is the restriction of $\G$ to $\d G$. If $M\ss X$ is finite and $M\ss\G(V)$, then we say that the network $\G$ {\em joins the set $M$}. The vertices of graphs and networks which are not boundary are said to be {\em interior}. The value
 $$
\smt(M)=\inf\{d(\G)\mid\text{$\G$ is a network joining $M$}\}
 $$
is called the {\em length of shortest network}. Notice that the network $\G$ which joins $M$ and satisfies $d(\G)=\smt(M)$ may not exist, see~\cite{ITLup}
and~\cite{Borod} for non-trivial examples. If such a network exists, it is called a {\em shortest network joining $M$}, or {\em for $M$}. One variant of the
Steiner problem is to describe the shortest networks joining finite subsets of pseudo-metric
spaces.%
 \footnote{%
The denotation $\smt$ is an abbreviation of ``Steiner Minimal Tree'' which is a synonym for the shortest network whose edges are non-degenerate and, thus, it must be a tree.
 }

Now let us define minimal parametric networks in a pseudo-metric space $\cX=(X,d)$. Let $G=(V,E)$ be a connected graph with some boundary $\d G$, and let $\v\:\d G\to X$ be a mapping. By $[G,\v]$ we denote the set of all networks $\G\:V\to X$ of the type $G$ such that $\d\G=\v$. We put
 $$
\mpn(G,\v)=\inf_{\G\in[G,\v]}d(\G)
 $$
and the value obtained we call by the {\em length of minimal parametric network}. If there exists a network $\G\in[G,\v]$ such that $d(\G)=\mpn(G,\v)$, then $\G$ is called a {\em minimal parametric network of the type $G$ with the boundary $\v$}.

 \begin {ass}
Let $\cX=(X,d)$ be an arbitrary pseudo-metric space and $M$ be a finite subset of $X$. Then
 $$
\smt(M)=\inf\{\mpn(G,\v)\mid\v(\d G)=M\}.
 $$
 \end {ass}

Thus, the problem of calculating the length of the shortest network is reduced to investigation of minimal parametric networks.

Let $\cM=(M,\r)$ be a finite pseudo-metric space and $G=(V,E)$ be an arbitrary connected graph joining $M$. In this case we always assume that the boundary of such graph $G$ is fixed and equal $M$. By $\Om(\cM,G)$ we denote the set of all weight functions $\om\:E\to\R$ such that
$(G,\om)$ is a filling of the space $\cM$. We put
 $$
\mpf(\cM,G)=\inf_{\om\in\Om(\cM,G)}\om(G)
 $$
and the value obtained we call the {\em weight of minimal parametric filling of the type $G$ for the space $\cM$}. If there exists a weight function $\om\in\Om(\cM,G)$ which $\om(G)=\mpf(\cM,G)$ for, then $(G,\om)$ is called a {\em minimal parametric filling of the type $G$ for the space $\cM$}.

 \begin {ass}
Let $\cM=(M,\r)$ be a finite pseudo-metric space. Then
 $$
\mf(\cM)=\inf\bigl\{\mpf(\cM,G)\bigr\}.
 $$
 \end {ass}

 \section {Reduction to trees} \label{sec:reduction}
 \markright {\thesection.~Reduction to trees.}
In this paragraph we show that investigating minimal fillings and shortest networks, one can restrict himself by trees with special boundaries. Also, it is sufficient to consider not all but some graphs, rather simple ones, in the study of minimal parametric fillings and networks. To start with, we give definitions of fillings and networks in the most general terms of multigraphs.

 \subsection {Multigraphs, fillings, and networks}
Recall that {\em multigraph\/} is a triple $G=(V,E,I)$, where $V$ and $E$ are finite sets which are called the sets of {\em vertices\/} and {\em edges\/} of $G$, respectively, and $I$ is a mapping from $E$ to the set of $1$- or $2$-element subsets of $V$. The mapping $I$ is called the {\em incidence}. If $I(e)=\{v\}$, then $e$ is called a {\em loop} incident to $v$. For simplicity reasons we sometimes denote the $1$-element set $\{v\}$ by  $\{v,v\}$. Moreover, we sometimes denote an edge $e\in E$ such that $I(e)=\{x,y\}$ by $xy$.

Further, if $I^{-1}(x)\ne\0$, then the number of elements in $I^{-1}(x)$ is called the {\em multiplicity\/} of each of the {\em edges\/} $e\in I^{-1}(x)$. Each edge with multiplicity more than $1$ is called {\em multiple}. A graph without loops and multiple edges is called {\em simple}. For such graphs the image of the mapping $I$ is a set of $2$-element subsets of $V$ (loops are absent) and $I$ is one-to-one with its image (multiple edges are absent), thus in this case $E$ may be identified with $I(E)$, i.e., with a family of $2$-element subsets of $V$, and the mapping $I$ is usually omitted. Notice that in previous paragraphs we handled just the simple graphs, and called them by graphs simply. In the current paragraph we use the term {\em graph\/} for multigraphs as a rule.

The concepts of paths, cycles, connectivity, weight function, filling, network can be defined for multigraphs in the same manner as for simple graphs. Let us note that in any network all its loops have zero lengths, and all multiple edges joining the same pair of vertices have the same lengths.

The next construction will be useful in what follows. Let $G=(V,E,I)$ be an arbitrary graph, $X$ be some set, and $f\:V\to X$ be a mapping. We use $f$ to denote the extension of $f$ to the set of all subsets of $V$ also. Define the {\em image $f(G)$ of the graph $G$} as $f(G)=\bigl(f(V),E,f\c I\bigr)$. If $M$ is the boundary of $G$, then the {\em boundary of graph $f(G)$} is the set $f(M)\ss X$.

Now we give a list of some simplest properties of minimal fillings and networks considering in such generality. These properties will show that such generality is redundant. We start from minimal parametric fillings. In what follow we call the edges of zero length by {\em degenerate\/} ones, and the edges of non-zero length by {\em non-degenerate\/} ones. A weighed graph without degenerate edges is said to be {\em non-degenerate}.

 \begin {ass}\label{ass:false1}
Let $\cG$ be a minimal parametric filling of the type $G$ of some pseudo-metric space $\cM$. Then
 \begin {description}
  \item[(1)] all interior vertices of degree $1$ of the graph $G$ are incident to degenerate edges\/\rom;
  \item[(2)] if one removes an interior vertex of degree $2$ incident to different edges from the graph $G$ and afterwards glues this two edges to one and assigns to the edge obtained the sum of weights of its generating edges, then the resulting weighted graph is a minimal parametric filling of $\cM$ as well\rom;
  \item[(3)] each loop of the graph $G$ is degenerate\/\rom;
  \item[(4)] multiple edges of $G$ have the same weights\/\rom;
  \item[(5)] the weight of each edge of $G$ being a part of a cycle does not exceed the sum of weights of the remaining edges of the cycle.
 \end {description}
 \end {ass}

 \begin {proof}
(1) Changing the weights of edges incident to vertices of degree $1$, one preserves the property of the weighted graph to be a filling, thus, by minimality reasons, such edges must have zero weights.

(2) Gluing the two edges incident to an interior vertex of degree $2$ to one edge of the total weight does not change as the property of graph to be a filling, so as the weight of the graph, thus the obtained graph must be a minimal parametric filling as well.

(3) The proof is similar with the one of (1).

(4) If multiple edges have different weights, then decreasing the biggest weight slightly, we preserve the property of the graph to be a filling, but decrease the total weight of the graph.

(5) Otherwise, one can decrease the weight of the edge preserving $G$ to be a filling. The proof is complete.
 \end {proof}

Now let us pass to minimal fillings. Since each of them is a minimal parametric filling, it has all the properties from Assertion~\ref{ass:false1}. However, some of the properties can be reinforced.

 \begin {ass}\label{ass:min_filling_cycles_mult_edges}
In any minimal filling
 \begin {description}
  \item[(1)] each multiple edge is degenerate\/\rom;
  \item[(2)] each edge contained in a cycle is degenerate.
 \end {description}
 \end {ass}

 \begin {proof}
From each family of non-degenerate multiple edges one can remove an edge and the resulting graph remains a filling but with less weight. Thus, in each minimal filling all multiple edges must be degenerate. The same arguments can be applied to non-degenerate edges in cycles. The proof is complete.
 \end {proof}

Similarly to the case of fillings, we call the zero-length edges of a network by {\em degenerate\/} ones, and all the remaining edges by {\em non-degenerate}. A network which has no degenerate edges is said to be {\em non-degenerate}.

 \begin {ass}\label{ass:false2}
Let $\G$ be a minimal parametric network of the type $G=(V,E,I)$ in a pseudo-metric space $(X,d)$. Then
 \begin {description}
  \item[(1)] all vertices of degree $1$ of the network $\G$ are incident to degenerate edges\/\rom;
  \item[(2)] each interior vertices of degree $2$ of the network $\G$ incident to different edges is placed between the vertices adjacent to it, i.e., if $v$, $vw_1$, and $vw_2$ are the corresponding vertices and edges of the parametric graph $G$, then
 $$
d\bigl(\G(w_1),\G(v)\bigr)+d\bigl(\G(v),\G(w_2)\bigr)=
d\bigl(\G(w_1),\G(w_2)\bigr),
 $$
thus, if one removes $v$ from $G$ and changes the edges $vw_i$ to the edge $w_1w_2$, then the network obtained by the restriction of $\G$ to the set $V\sm\{v\}$ of the reconstructed graph $G$ is a minimal parametric network with the same boundary\rom;
  \item[(3)] each loop of the network $\G$ \(as the one of any other network\/\) is degenerate\/\rom;
  \item[(4)] multiple edges of the network $\G$ \(as the ones of any other network\/\) have the same length.
 \end {description}
 \end {ass}

 \begin {proof}
(1) Suppose on the contrary that there exists an interior vertex $\G\:v\to X$ of degree $1$ incident to an edge $\G\:vw\to X$ of a positive length. Then we change the mapping $\G$ at the vertex $v$ as $\G(v)=\G(w)$ and obtain a network with the same type and boundary but with less length, a contradiction.

(2) Once more, suppose on the contrary that there exists an interior vertex $\G\:v\to X$ of degree $2$ such that for its adjacent vertices $\G\:w_1\to X$ and $\G\:w_2\to X$ we have
 $$
d\bigl(\G(w_1),\G(v)\bigr)+d\bigl(\G(v),\G(w_2)\bigr)
\ne d\bigl(\G(w_1),\G(w_2)\bigr).
 $$
Then, by triangle inequality, we obtain
 $$
d\bigl(\G(w_1),\G(v)\bigr)+d\bigl(\G(v),\G(w_2)\bigr)
>d\bigl(\G(w_1),\G(w_2)\bigr).
 $$
Thus, if we change the mapping $\G$ at the vertex $v$ in such a way that $v$ goes to $\G(w_1)$, then we obtain a parametric network with the same boundary but with less length, a contradiction.

Items (3) and (4) are evident. The proof is complete.
 \end {proof}

Now we pass to shortest networks. Since each of them is a minimal parametric network, it has all the properties from Assertion~\ref{ass:false2}. However, some of them can be reinforced.

 \begin {ass}\label{ass:smt}
In each shortest network
 \begin {description}
  \item[(1)] any multiple edge is degenerate\/\rom;
  \item[(2)] each edge contained in a cycle is degenerate.
 \end {description}
 \end {ass}

 \begin {proof}
From each family of non-degenerate multiple edges one can remove an edge and the resulting parametric graph is connected and join points from $M$, but the resulting network will be shorter. Thus, in each shortest network all multiple edges must be degenerate. The same arguments can be applied to non-degenerate edges in cycles. The proof is complete.
 \end {proof}

From Assertions~\ref{ass:false1} and~\ref{ass:false2} we can see that removing interior vertices of degree $2$ incident to different edges does not violate the minimality. It is easy to define the inverse operation which is called {\em edge subdivision}. Namely, the operation consists in adding an interior vertex of degree $2$ inside an edge. Clearly, the operation preserves minimality as well. All that motivates the term {\em false \(redundant\/\)} used to call the interior vertices of degree $2$ incident to different edges. As a rule, we suppose that the graphs and networks in consideration do not contain false vertices.

The next step on the way of simplification of the networks to work with is to forbid the degenerate edges described in the previous Assertions. To do that, we define factorization of a graph $G$ by a family of its edges.

 \subsection {Factorization and splitting}
Let $G=(V,E,I)$ be an arbitrary graph with a boundary $M$, and $F$ be a subset of $E$. By $G_i=(V_i,F_i,I_i)$ we denote the connected components of the graph $(V,F,I)$ and put $M_i=M\cap V_i$, $V_F=\{V_i\}$, $M_F=\{M_i\}$. Notice that $V_F$ и $M_F$ are partitions of the sets $V$ and $M$, respectively. Let $\pi_F\:V\to V_F$ be a natural projection which is uniquely defined by the condition $v\in\pi_F(v)$ for each $v\in V$. By $\pi_F$ we denote also the restriction of this projection to $M$ and the extension of this projection to the set of all subsets of $V$. We put $E_F=E\sm F$, and for each $e\in E_F$ we define $I_F$ as follows: $I_F(e)=\pi_F\bigl(I(e)\bigr)$.  The graph  $(V_F,E_F,I_F)$ with the boundary $M_F$ is called the {\em factor\/} of the graph $G$ by the edges family $F$ and is denoted by $G_F$.

Now, we transfer the concept of factorization to fillings and networks. To do that, let us consider an arbitrary pseudo-metric space $\cM=(M,\r)$, and let $P=\{M_i\}$ be a partition of the set $M$. The partition $P$ is called {\em metric}, if the restriction of the distance $\r$ to each $M_i$ vanishes. It is easy to see that in this case for any $i$ and $j$ the distance $\r(x,y)$ between $x\in M_i$ and $y\in M_j$ does not depend on the choice of $x$ and $y$. Moreover, the function $\r_P(M_i,M_j)$ which is equal to $\r(x,y)$ for some $x\in M_i$ and $y\in M_j$ defines a pseudo-metric on $P$. We say that the function $\r_P$  {\em is induced\/} by the metric partition $P$. We call the pseudo-metric space $\cM_P=(P,\r_P)$ by the {\em factor\/} of the space $\cM$ by the metric partition $P$. Notice that for metric partition $P$ the projection $\pi_P\:M\to P$ which is uniquely defined by the condition $x\in\pi_P(x)$ preserves the distances, i.e., $\r(x,y)=\r_P\bigl(\pi_P(x),\pi_P(y)\bigr)$ for any $x$ and $y$ from $M$.

Let $\cG=(G,\om)$ be a filling of a type $G=(V,E,I)$ of a pseudo-metric space $\cM=(M,\r)$, and suppose that $F\ss E$ consists of degenerate edges, probably, not from all ones.  As above, let $G_F=(V_F,E_F,I_F)$ be the factor of the graph $G$ by the family $F$, and $M_F$ be its boundary. By $\om_F$ we denote the restriction of $\om$ to $E_F=E\sm F$, and by $\cG_F$ the weighted graph $(G_F,\om_F)$. Since all edges of each connected component $G_i=(V_i,F_i,I_i)$ of the graph $(V,F,I)$ are degenerate, the restriction of $d_\om$ to each $V_i\in V_F$ vanishes, thus, since $\cG$ is a filling, the restriction of $\r$ to each $M_i\in M_F$ vanishes as well. Thus, the partition $M_F$ of the set $M$ is metric, therefore the factor $\cM_F=(M_F,\r_F)$ by the partition $M_F$ is defined. It is easy to verify that $\cG_F$ is a filling of $\cM_F$. The filling $\cG_F$ is called the {\em factor\/} of the filling $\cG$ by the family $F$ of degenerate edges.

Now, let $\G\:V\to X$ be a network of a type $G=(V,E,I)$ with some boundary $\v\:N\to X$, joining the set $M=\v(N)$ in a pseudo-metric space $\cX=(X,d)$, and let $F\ss E$ consist of some degenerate edges. Let $G_F$ be a graph with the boundary $N_F$ such that $G_F$ is the factor of $G$ by the family $F$. Consider the graph $H=\G(V,F,I)$ and let $H_k=(W_k,F_k,J_k)$ be the connected components of this graph. Construct a partition $W$ of the set $X$ by extending the family $\{W_k\}$ with pointwise subsets of $X$, whose elements do not belong to any $W_k$. It is easy to see that $W$ is a metric partition of the space $\cX$. Let $P$ be an arbitrary partition of the space $\cX$ such that $P$ is an enlargement of $W$, i.e., $W$ is a subdivision of $P$. In this case we say that the family $F$ is {\em compatible\/} with the metric partition $P$ of the space $\cX$, and we call the pair $D=(F,P)$ by a {\em compatible degenetarion\/} of the network $\G$. Let $\cX_P$ be the factor of $\cX$ by the partition $P$, and $\pi_P\:\cX\to\cX_P$ be the corresponding projection. Since $P$ is an enlargement of $W$, then the set $\pi_P\bigl(\G(V_i)\bigr)$ is pointwise for each $V_i\in V_F$. Now, we define the network $\G_D\:V_F\to\cX_P$ with the boundary $\v_D\:N_F\to\cX_P$, parameterized by the graph $G_F=(V_F,E_F,I_F)$. To do that, we put $\G_D(V_i)=\pi_P\bigl(\G(V_i)\bigr)$ and define $\v_D$ to be the restriction of $\G_D$ to $N_F$. The obtained network $\G_D$ joins the set $M_P=\pi_P(M)=\v_D(N_F)$.

Notice, that for fillings of {\sl metric\/} spaces and networks in metric spaces the definitions of factorization are much more simple because we do not need to factorize the spaces. Let us show the corresponding constructions.

Let $\cG=(G,\om)$ be a filling of a type $G=(V,E,I)$ of a metric space $\cM=(M,\r)$, and $F\ss E$ consist of some degenerate edges. Since the space $\cM$ is a metric one, the partition $M_F$ consists of pointwise subsets of $M$ and, thus, $M_F$ is naturally identified with $M$. Therefore, in the case of metric spaces the weighted graph $\cG_F$ constructed above is a filling of the initial space $\cM$.

Now, let $\G\:V\to X$ be a network of a type $G=(V,E,I)$ with a boundary $\v\:N\to X$ in a metric space $\cX=(X,d)$, and $F\ss E$ consists of some degenerate edges. Notice that in this case an edge $e\in E$ is degenerate iff $\G\bigl(I(e)\bigr)$ is a pointwise subset of $X$, thus the sets $\G(V_i)$ are pointwise for each $V_i\in V_F$. Moreover, the space $\cM$ has the unique metric partition, namely, the partition onto pointwise subsets. Therefore, in this case the factor-network $\G_D\:V_F\to X$ with the boundary $\v_D\:N_F\to X$ can be defined as follows: $\G_D(V_i)=\G(V_i)$ and $\v_D$ is the restriction of $\G_D$ to $N_F\ss V_F$.

The inverse operation to factorization is called splitting. Namely, if a graph, a weighted graph, a network, and the corresponding pseudo-metric spaces are obtained by factorization from the corresponding objects, then we say that the initial object are obtained from the terminal one by {\em splitting}.

 \subsection {Base}
Let $G=(V,E,I)$ be a graph, $\cX=(X,d)$ be a pseudo-metric space, $\cG=(G,\om)$ be a weighted graph, $\G\:V\to X$ be a network in $\cX$. By $\degen(\cG)$ and $\degen(\G)$ we denote the sets of all degenerate edges of the weighted graph $\cG$ and the network $\G$, respectively. The graph $\cG_F$ and the network $\G_F$ for $F=\degen(\cG)$ and $F=\degen(\G)$, respectively, are called the {\em subbases\/} of the graph $\cG$ and network $\G$, respectively.  Subbase from which the false vertices removed is called the {\em base}. Notice that the transformation to subbase preserves the weight of the graph and the length of the network, respectively, but the transformation to base preserves the weight only. The length may decrease.

For a pseudo-metric space $\cX$ by $\degen(\cX)$ we denote the unique metric partition from which one can obtain all other metric partitions by means of iterative use of subdivision operation. The elements of the partition $\degen(\cX)$ are called the {\em degenerate components\/} of the space $\cX$. The {\em base\/} of the space $\cX$ is the factor-space $\cX_P$, where $P=\degen(\cX)$.

 \begin {ass}\label{ass:false_gen1}
Let $\cM$ be a pseudo-metric space, and $\cG$ be some its minimal\/ \(parametric\/\) filling.
 \begin {description}
  \item[(1)] If $F$ is an arbitrary family of degenerate edges from $\cG$, then the weighted graph $\cG_F$ is a minimal filling\/ \(minimal parametric filling, respectively\/\) of the space $\cM_F$.
  \item[(2)] If $F=\degen(\cG)$, then the base $\cG_F$ of the minimal\/ \(parametric\/\) filling $\cG$ of the space $\cM$ is minimal\/ \(parametric\/\) filling of the space $\cM_F$ and it is a non-degenerate graph without loops with the boundary containing all the vertices of degree $1$ and $2$. Besides, if $\cG$ is a minimal filling, then its base $\cG_F$ is a non-degenerate tree and the space $\cM_F$ coincides with the base of the space $\cM$, in particular,  $\cM_F$ is a metric space.
 \end {description}
 \end {ass}

 \begin {proof}
{\bf (1)} From definition of the graph $\cG_F$ it follows that the graph is a filling of the space $\cM_F$ and its weight equals to the one of the initial graph $\cG$, therefore $\cG_F$ is minimal (parametric) filling of $\cM_F$.

{\bf (2)} Notice that removing of the false vertices preserves the minimality, and the presence in the base of non-degenerate loops or non-degenerate edges incident to interior vertices of degree $1$ contradicts to Assertion~\ref{ass:false1}.

By Assertion~\ref{ass:min_filling_cycles_mult_edges}, all the edges belonging to cycles, and all multiple edges, are degenerate, thus the minimal filling $\cG_F$ is a non-degenerate tree. It remains to show that $\cM_F$ in this case coincides with the base of the space $\cM$.  On the contrary, suppose that the space $\cM_F$ contains different points $p$ and $q$ for which $\r_F(p,q)=0$, and let $\g$ be a path in $\cG_F$ joining $p$ and $q$. The path $\g$ contains an edge $e$ of some weight $\om>0$, since otherwise the points $p$ and $q$ are identified in $M_F$. If $e\ne pq$, then we reconstruct the tree $\cG_F$ by removing the edge $e$ and adding the one $pq$ of the weight $\om$. The reconstructed graph is a minimal filling as well.

 \begin {lem}\label{lem:exact_edge}
Let $\cG=(G,\om)$ be a minimal parametric filling of a pseudo-metric space $(M,\r)$, and $e$ be an edge from $G$ joining some points $p$ and $q$ from $M$. Then $\r(p,q)=\om(e)$.
 \end {lem}

 \begin {proof}
Suppose on the contrary that $\r(p,q)\ne\om(e)$. Since $\cG$ is a filling, then we have $\r(p,q)<\om(e)$. Let $\g$ be an arbitrary path in $\cG$ containing $e$ and joining some points from $M$, say $p'$ and $q'$. To be definite, suppose that along the path $\g$ the points in consideration are placed in the order $p'$, $p$, $q$, $q'$. Then
 $$
\r(p',q')\le\r(p',p)+\r(p,q)+\r(q,q')<\r(p',p)+\om(e)+\r(q,q')\le\om(\g),
 $$
therefore there exists sufficiently small positive $\e$ such that the tree obtained from $\cG$ by decreasing the weight of the edge $e$ by $\e$ remains a filling, that contradicts to minimality of the filling $\cG$. The proof is complete.
 \end {proof}

Lemma~\ref{lem:exact_edge} implies that the weight of the edge $e$ vanishes, a contradiction which completes the proof of Assertion.
 \end {proof}

 \begin {rk}
For the base of minimal parametric filling the corresponding factor of pseudo-metric space need not, generally speaking, be a metric space. Figure~\ref{fig:non_metric}, left, shows an example of pseudo-metric space $\cM$ generated by weighted graph. The space $\cM$ is not metric. In the center it is depicted a minimal parametric filling $\cG$ of the space $\cM$. Thus, the base $\cG_F$ of this filling coincides with the filling, but the space $\cM_F$ is not a metric one. The minimal filling $\cH$ of this space has another type and is depicted in Figure~\ref{fig:non_metric}, right. Notice that for the base $\cH_R$ of this filling, the space $\cM_R$ does be a metric one, what is concordant to Assertion~\ref{ass:false_gen1}.

\Pic{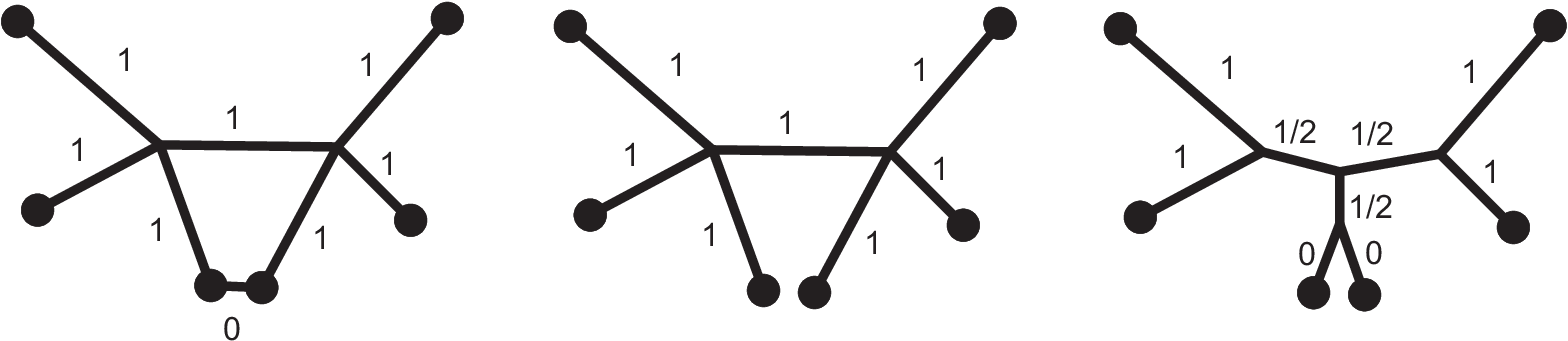}{Example of pseudo-metric space which is not a metric one, but which is the boundary of the base of a minimal parametric filling.}{fig:non_metric}{100}
 \end {rk}

Similarly to Assertion~\ref{ass:false_gen1}, one can prove its following variant for the case of the shortest networks.

 \begin {ass}\label{ass:false_gen2}
Let $\cX=(X,d)$ be a pseudo-metric space and $\G$ a shortest or minimal parametric network in $\cX$ joining $M\ss X$.
 \begin {description}
  \item[(1)] If $D=(F,P)$ is an arbitrary compatible degeneration of the network $\G$, then $\G_D$ is a shortest\/ \(minimal parametric\/\) network in the space $\cX_P$ joining the set $M_P$.
  \item[(2)] If $F=\degen(\G)$ and $D=(F,P)$ is a compatible degeneration of the network $\G$, then the base $\Y$ of the network $\G$ is a shortest\/ \(minimal parametric\/\) non-degenerate network in the metric space $\cX_P$ joining the set $M_P$. Also, the network $\Y$ does not contain loops, and int boundary contains all the vertices of degree $1$ and $2$. Besides, the base $\Y$ of the shortest network is a non-degenerate tree, and the restriction of the pseudo-metric on $M_P\ss X_P$ is a metric.
 \end {description}
 \end {ass}

Notice that the base of each minimal filling and of each shortest network can be obtained by factorization procedure from some tree $T$ which has only vertices of degree $1$ and $3$, and whose boundary coincides with the boundary of the base and consists just of all vertices from $T$ of degree $1$. We call such trees {\em binary}. Thus, we have obtained the following result.

 \begin {ass}
The base of a minimal filling or of a shortest network can be splitted to a binary tree which is, respectively, a minimal filling or a shortest network with the same boundary.
 \end {ass}

To summarize, we formulate the following result.

 \begin {thm}\label{th:trees_fill}
For any finite pseudo-metric space $\cM$, there exists a minimal filling which is a binary tree. If the space $\cM$ is a metric one, then it has as well a minimal filling which is a non-degenerate tree. Conversely, if for $\cM$  there exists a non-degenerate minimal filling, then $\cM$ is a metric space.

For any finite subset $M$ of a pseudo-metric space, there exists a shortest binary tree joining $M$. If the distances between different points from $M$ are non-zero, then there exists a shortest non-degenerate tree joining $M$. Conversely, if there exists a shortest non-degenerate tree  for $M$, then the distances between different points from $M$ are non-zero.
 \end {thm}

 \subsection*{Agreement}
Assertions~\ref{ass:false_gen1} and~\ref{ass:false_gen2} imply that to investigate minimal filling and shortest networks one may consider only the trees whose vertices of degree $1$ and $2$ belong to the boundary. {\bf In what follows, we always assume that this condition holds, providing the opposite is not declared}.

 \section {Minimal realization}\label{sec:realization}
 \markright {\thesection.~Minimal realization.}
In this Section we show that the problem on minimal filling can be reduced to Steiner problem in special metric spaces and for special boundaries.

Consider a finite set $M=\{p_1,\ldots,p_n\}$, and let $\cM=(M,\r)$ be a pseudo-metric space. We put $\r_{ij}=\r(p_i,p_j)$. By $\ell_\infty^n$  we denote the $n$-dimensional arithmetic space with the norm
 $$
\bigl\|(v^1,\ldots,v^n)\bigr\|_\infty=\max\bigl\{|v^1|,\dots,|v^n|\bigr\},
 $$
and by $\r_\infty$ the metric on $\ell_\infty^n$ generated by $\|\cdot\|_\infty$, i.e., $\r_\infty(v,w)=\|w-v\|_\infty$. Let us define a mapping $\v_\cM\:M\to\ell_\infty^n$ by the formula
 $$
\v_\cM(p_i)=\bp_i=(\r_{i1},\ldots,\r_{in}).
 $$

 \begin {ass}\label{ass:isom_embedding_ellinfty}
The mapping $\v_\cM$ is an isometry with its image.
 \end {ass}

 \begin {proof}
Indeed,
 $$
\r_\infty(\bp_i,\bp_j)=
\max\{|\r_{i1}-\r_{j1}|,\ldots,|\r_{ij}-\r_{jj}|,\ldots,|\r_{in}-\r_{jn}|\}
=\r_{ij}=\r(p_i,p_j),
 $$
because $|\r_{ik}-\r_{jk}|\le \r_{ij}$ и $\r_{jj}=0$. The  proof is complete.
 \end {proof}

The mapping $\v_\cM$ is called by the {\em Kuratowski isometry}.

Let $\cG=(G,\om)$ be a filling of a space $\cM=(M,\r)$, where $G=(V,E)$, and $d_\om$ be a pseudo-metric on $V$ generated by the weight function $\om$.  Consider the multigraph $\bG=(V,\bE,I)$, where $\bE$ is the disjoint union of the edge set $E$ of the tree $G$ and the edge set $E_M$ of the complete graph $K_n(M)$ on $M$. The incidence mapping $I$ is the union of incidence mapping of the simple graphs $G$ and $K_n(M)$. Let $\bom$ be the weight function on $\bE$ coinciding with $\om$ on $E$ and with pseudo-metric $\r$ on $E_M$. Notice that the edges from $\bG$ joining points from $M$ may have multiplicity $2$, but all the other edges are of the  multiplicity $1$. By $d_\bom$  we denote the pseudo-metric on $V$ generated by $\bom$.

 \begin {ass}\label{ass:ext-metric-to-filling}
We have
 \begin {description}
  \item[(1)] $d_\bom(p_i,p_j)=\r(p_i,p_j)\le d_\om(p_i,p_j)$\rom;
  \item[(2)] for any $v$ and $w$ from $V\sm M$ we have
 $$
d_\bom(v,w)=\min_{\a,\b}\bigl\{d_\om(v,p_\a)+\r(p_\a,p_\b)
+d_\om(p_\b,w),d_\om(v,w)\bigr\}\le d_\om(v,w),
 $$
in particular, if $e=vw\in E$, then $d_\bom(v,w)\le\om(e)$\rom;
  \item[(3)] for any $v\in V\sm M$ we have
 $$
d_\bom(v,p_i)=\min_{\a}\bigl\{d_\om(v,p_\a)+\r(p_\a,p_i)\bigr\}\le
d_\om(v,p_i).
 $$
 \end {description}
 \end {ass}

 \begin {proof}
In each of the three cases in consideration, by $\g$ we denote a path in $\bG$ joining the vertices from the assumption. Partition the path $\g$ into consecutive paths  $\g_i$, $i=1,\ldots,k$, where $\g_i$ is a maximal sub-path of $\g$ consisting of either just the edges from $E_M$ (the {\em first type\/} path), or of the edges from $E$ (the {\em second type\/} path). Notice that for each $i=2,\ldots,k-1$ the path $\g_i$ joins some vertices $p_\a$ and $p_\b$ from $M$. If such $\g_i$ belongs to the first type, then $\bom(\g_i)\ge\r(p_\a,p_\b)$ by the triangle inequality. If the path $\g_i$ belongs to the second type, then $\bom(\g_i)\ge\r(p_\a,p_\b)$ because $(G,\om)$ is a filling. Thus, if $\g'$ denotes the union of the paths $\g_i$ over all $i=2,\ldots,k-1$, and $p_\a$, $p_\b$ are its ending vertices, then $\bom(\g')\ge\r(p_\a,p_\b)$ by the triangle inequality.

Now, consider the first item of Assertion. In this case $\g_1$ and $\g_k$ join as well the vertices from $M$, therefore by similar reasons we get that $\bom(\g)\ge\r(p_i,p_j)$, and thus $d_\bom(p_i,p_j)\ge\r(p_i,p_j)$. Notice that in this case, among such $\g$ there exists a path consisting of the edge $p_ip_j\in E_M$ whose weight equals $\r(p_i,p_j)$, therefore $d_\bom(p_i,p_j)=\r(p_i,p_j)$.

In the second Item we have either $\g=\g_1$ and, thus, $\bom(\g)=\om(\g_1)\ge d_\om(v,w)$, or $k\ge 2$, $\bom(\g_1)=\om(\g_1)\ge d_\om(v,p_\a)$, and $\bom(\g_k)=\om(\g_k)\ge d_\om(p_\b,w)$. Besides, the previous results imply that $d_\bom(p_\a,p_\b)=\r(p_\a,p_\b)$, thus
 $$
d_\bom(v,w)=\min_\g\bigl\{\bom(\g)\bigr\}\ge
\min_{\a,\b}\bigl\{d_\om(v,p_\a)+\r(p_\a,p_\b)+d_\om(p_\b,w),
d_\om(v,w)\bigr\}.
 $$
If
 $
\min_{\a,\b}\bigl\{d_\om(v,p_\a)+\r(p_\a,p_\b)+d_\om(p_\b,w)\bigr\}
\ge d_\om(v,w)
$, then we choose $\g$ to be the path in $G$, which $\om(\g)=d_\om(v,w)$ for. If the inverse inequality holds, then we put $\a$ and $\b$ be those indices which
the minimum is attained at, and $\g$ be the path composed from $\g_1$, $\g_2$ and $\g_3$, where $\g_1$ and $\g_3$ are the paths in $G$ for which $\om(\g_1)=d_\om(v,p_\a)$, $\om(\g_3)=d_\om(p_\b,w)$, and the path $\g_2$ consists of just one edge $p_\a p_\b\in E_M$. In each of these cases the inequality turns to the equality sought for.

The third Item can be proved similarly to the second one. The proof is complete.
 \end {proof}

We define the network $\G_\cG\:V\to\ell_\infty^n$ of the type $G$ as follows:
 $$
\G_\cG(v)=\bigl(d_{\bom}(v,p_1),\ldots,d_{\bom}(v,p_n)\bigr).
 $$
This network is called by the {\em Kuratowski network for the filling $\cG$}.

Assertion~\ref{ass:ext-metric-to-filling} implies the following result immediately.

 \begin {ass}\label{ass:extend}
We have $\d\G_\cG=\v_\cM$.
 \end {ass}

A mapping of a pseudo-metric spaces is called {\em non-stretching} if it does not increase the distances between points. From Assertion~\ref{ass:ext-metric-to-filling} we get one more result.

 \begin {ass}\label{ass:non-stratch}
The mapping $\G_\cG$ acting from $(V,d_\om)$ to $\ell_\infty^n$ is non-stretching.
 \end {ass}

 \begin {proof}
Let $v$ and $w$ be arbitrary vertices from $V$, then
 \begin {multline*}
\r_\infty\bigl(\G_\cG(v),\G_\cG(w)\bigr)=\\ =
\max
\bigl\{
|d_\bom(v,p_1)-d_\bom(w,p_1)|,\ldots,|d_\bom(v,p_n)-d_\bom(w,p_n)|
\bigr\}\le\\
\le d_\bom(v,w)\le d_\om(v,w),
 \end {multline*}
where the first inequality holds by the triangle rule, and the last inequality --- by Assertion~\ref{ass:ext-metric-to-filling}. The proof is complete.
 \end {proof}

 \begin {rk}
An analogue of Assertion~\ref{ass:non-stratch} takes place in much more general assumptions, see~\cite{SIvanovDis}.
 \end {rk}

For any network $\G$ in a metric space $(X,d)$ by $\om_\G$ we denote the {\em weight function on $G$ induced by the network $\G$},
i.e., $\om_\G(vw)=d\bigl(\G(v),\G(w)\bigr)$.

 \begin {cor}\label{cor:induced_from_ell}
Let $\cG=(G,\om)$ be a minimal parametric filling of a metric space $(M,\r)$ and $\G=\G_\cG$ the corresponding Kuratowski network. Then $\om=\om_\G$.
 \end {cor}

 \begin {proof}
Since $\G$ maps isometrically the space $(M,\r)$ to $\ell_\infty^n$, the weighted graph $(G,\om_\G)$ is a filling of the space $(M,\r)$. By Assertion~\ref{ass:non-stratch}, for each edge $e\in E$ we have $\om(e)\ge\om_\G(e)$. Therefore, by minimality of the filling $(G,\om)$, all these inequalities are satisfied in the form of equalities. The proof is complete.
 \end {proof}

Let $\G$ be a network in a pseudo-metric space $\cX$, let $G$ be its parameterizing graph, and $\cH=(H,\om)$ be a weighted graph. We say that {\em $\G$ and $\cH$ are isometric}, if there exists an isomorphism of the weighted graphs $\cH$ and $\cG=(G,\om_\G)$.

Corollary~\ref{cor:induced_from_ell} and the existence of minimal parametric and shortest networks in a finite-dimensional normed space~\cite{ITCRC} imply the following result.

 \begin {cor}\label{cor:Kurat_mpf}
Let $\cM=(M,\r)$ be a pseudo-metric space  consisting of $n$ points, and $\v_\cM\:M\to\ell_\infty^n$ be the Kuratowski isometry. For any graph $G$ joining $M$ there exists a minimal parametric filling of the type $G$ of the space $\cM$. Each minimal parametric filling of the type $G$ of the space $\cM$ is isometric to the corresponding Kuratowski network, which is, in this case, a minimal parametric network of the type $G$ with the boundary $\v_\cM$. Conversely, each minimal parametric network of the type $G$ on $\v_\cM(M)$ is isometric to some minimal parametric filling of the type $G$ of the space $\cM$.
 \end {cor}

 \begin {cor}\label{cor:Kur_mf}
Let $\cM=(M,\r)$ be a pseudo-metric space  consisting of $n$ points, and $\v_\cM\:M\to\ell_\infty^n$ be the Kuratowski isometry. Then there exists a minimal filling $\cG$  for $\cM$, and the corresponding Kuratowski network $\G_\cG$ is a shortest network in the space $\ell_\infty^n$ joining the set $\v_\cM(M)$. Conversely each shortest network on $\v_\cM(M)$ is isometric to some minimal filling of the space $\cM$.
 \end {cor}

 \section {Reformulation of the Problem} \label{sec:reformul}
 \markright {\thesection.~Reformulation of the Problem.}
Let $\cM=(M,\r)$ be a finite pseudo-metric space, $\cX=(X,d)$ an arbitrary pseudo-metric space, and $f\:M\to X$ be an isometry with its image. For such mappings it is convenient to write $f\:\cM\to\cX$. By $\cI(\cM)$ we denote the set of all such $f\:\cM\to\cX$, provided $\cM$ is fixed, but $\cX$ is not.

 \begin {ass}\label{ass:FillingInTermsSTM}
We have $\mf(\cM)=\inf_{f\in\cI(\cM)}\smt\bigl(f(M)\bigr)$.
 \end {ass}

 \begin {proof}
Notice that each shortest tree on $f(M)$ is isometric to some filling of $\cM$, therefore
 $$
\mf(\cM)\le\inf_{f\in\cI(\cM)}\smt\bigl(f(M)\bigr).
 $$
If we take the mapping $\v\:\cM\to\ell_\infty^n$ defined above as $f$, then we obtain $\mf(\cM)=\smt\bigl(\v(M)\bigr)$. The proof is complete.
 \end {proof}

Let $\cM=(V,\r)$ be an arbitrary finite pseudo-metric space and $G=(V,E)$ be some graph. Such graph is called a {\em graph on the space $\cM$}. The {\em length $\r(e)$ of an edge $e=vw\in E$} is the value $\r(v,w)$, and the {\em length $\r(G)$ of the graph $G$} is the sum of lengths of all the graph edges. The value
 $$
\mst(\cM)=\min\{\r(G)\mid\text{$G$ is a tree on $\cM$}\}
 $$
is called the {\em length of minimal spanning tree\/} of the pseudo-metric space $\cM$, and the tree on $\cM$ for which $\r(G)=\mst(\cM)$ is called {\em minimal spanning tree\/} on the space $\cM$.

From the definitions one can immediately conclude the following result.

 \begin {ass}\label{ass:SMTinTermsMST}
Let $\cX=(X,d)$ be an arbitrary pseudo-metric space and $M$ be a finite subset of $X$, then
 $$
\smt(M,d)=
\inf\bigl\{\mst(N,d)\mid\text{$M\ss N\ss X$, $\#N<\infty$}\bigr\}.
 $$
 \end {ass}

From Assertions~\ref{ass:SMTinTermsMST} and~\ref{ass:FillingInTermsSTM} we get the following

 \begin {cor}\label{cor:mfInMST}
For any finite pseudo-metric space $\cM=(M,\r)$ the equality
 $$
\mf(\cM)=\inf_{f, N}\Bigl\{\mst(N,d)\mid
f\:\cM\to \cX,\ \cX=(X,d),\ f(M)\ss N\ss X
\Bigr\}
 $$
holds.
 \end {cor}

If $\cM$ and $\cN$ are finite pseudo-metric spaces such that an isometry $f\:\cM\to\cN$ with its image exists, then $\cN$ is called an {\em extension of $\cM$}. By $\cE(\cM)$ we denote the set of all extensions of the space $\cM$. It is easy to see that the spaces $(N,d)$ from Corollary~\ref{cor:mfInMST} run the whole set $\cE(\cM)$. This enables us to reformulate Corollary~\ref{cor:mfInMST} as follows.

 \begin {cor}\label{cor:secondDefOfMF}
The equality
 $$
\mf(\cM)=\inf_{\cN\in\cE(\cM)}\mst(\cN)
 $$
is valid.
 \end {cor}

Notice that the formula from Corollary~\ref{cor:secondDefOfMF} can be taken as an equivalent definition of the minimal filling weight.

Let $\cX=(X,\r)$ and $\cY=(Y,d)$ be pseudo-metric spaces. A mapping $f\:X\to Y$  is called {\em non-contracting}, if for any points $u$ and $v$ from $X$ the inequality $d\bigl(f(u),f(v)\bigr)\ge\r(u,v)$ holds. By $\cL(\cM)$ we denote the set of all finite pseudo-metric spaces $\cN$ such that there exists a non-contracting mapping of   a fixed pseudo-metric space $\cM$ into $\cN$.

 \begin {cor}\label{cor:thirdDefOfMF}
The equality
 $$
\mf(\cM)=\inf_{\cN\in\cL(\cM)}\mst(\cN)
 $$
is valid.
 \end {cor}

Corollary~\ref{cor:thirdDefOfMF} may be considered as the third equivalent definition of the minimal filling weight.

 \section {Minimal Parametric Fillings and Linear Programming} \label{sec:exist}
 \markright {\thesection.~Minimal Parametric Fillings and Linear Programming.}
Let $\cM=(M,\r)$ be a finite pseudo-metric space joined by a (connected) graph  $G=(V,E)$. As above, by $\Om(\cM,G)$ we denote the set consisting of all the weight functions  $\om\:E\to\R_+$ such that $\cG=(G,\om)$ is a filling of the space $\cM$, and by $\Om_m(\cM,G)$ we denote its subset consisting of the weight functions such that  $\cG$ is a minimal parametric filling of the space $\cM$.

 \begin {ass}\label{ass:opt_weight}
The set $\Om(\cM,G)$ is closed and convex in the linear space $\R^E$ of all the functions on $E$, and $\Om_m(\cM,G)\ss\Om(\cM,G)$ is a non-empty convex compact.
 \end {ass}

 \begin {proof}
Indeed, the condition that a weight function $\om$ defines a filling of the space  $\cM$ can be formulated as the following system of inequalities $\sum_{e\in\g}\om(e)\ge\r(p,q)$, where $\g$ is a path in $G$ joining vertices $p$ and $q$ from $M$, together with the inequalities
$\om(e)\ge0$, $e\in E$, which guarantee that the weights of edges are non-negative. Thus, the corresponding subset  $\Om(\cM,G)$ of the linear space
$\R^E$ is closed and convex.

Further, minimal parametric fillings correspond to minimums of the function $\sum_{e\in E}\om(e)$ on $\Om(\cM,G)$. Since this function is continuous and bounded from below on $\Om(\cM,G)$, the set $\Om_m(\cM,G)$ of minimal parametric fillings is non-empty and closed. It remains to notice that the level set of the function $\sum_{e\in E}\om(e)$ on the positive orthant is a simplex, and its level set on  $\Om(\cM,G)$ is the intersection of the simplex with the convex set $\Om(\cM,G)$, therefore $\Om_m(\cM,G)$ is convex and bounded. The proof is complete.
 \end {proof}

 \begin {rk}
The proof of Assertion~\ref{ass:opt_weight} shows that the problem of searching for a minimal parametric filling of a pseudo-metric space can be reduced to a linear programming.
 \end {rk}

 \section {Exact Paths, Exact Edges, Exact Vertices} \label{sec:exact}
 \markright {\thesection.~Exact Paths, Exact Edges, Exact Vertices.}
Let $\cM=(M,\r)$ be an arbitrary finite pseudo-metric space, and $\cG=(G,\om)$ be some its filling. A path joining boundary vertices of the tree  $\cG$ is called  {\em boundary}, and a boundary path whose weight is equal to the distance between its ending vertices is said to be  {\em exact}.

Let us start with the following useful technical remark. Let  $\cG=(G,\om)$ be a filling, $F=\degen(\cG)$ be the set of all degenerate edges from $\cG$, and  $\cG_F=(G_F,\om_F)$ be the base of this filling. Recall that by $\pi_F$ we denote the natural projection from the vertex set of the graph $G$ onto the vertex set of the graph $G_F$ extended to the subsets.

 \begin {ass}\label{ass:exact_project}
Under the above assumptions, a boundary path  $\g$ in $\cG$ is exact, if and only if the boundary path  $\pi_F(\g)$ in $\cG_F$ is exact.
 \end {ass}

 \subsection {General Fillings}

 \begin {ass}\label{ass:boud_paths}
If a boundary path in an arbitrary filling is contained in an exact path, then this path is also exact.
 \end {ass}

 \begin {proof}
Let $\dl$ be a boundary path which is contained in an exact path $\g$. By $x$ and $y$ we denote the ending vertices of the path $\g$, and by $z$ and $w$ we denote the ending vertices of the path $\dl$. Without loss of generality, assume that the vertices $x$, $z$, $w$, $y$ are consecutive along $\g$.  By $\g_x$ and $\g_y$ we denote the segments of the path $\g$ from $x$ to $z$ and from $w$ to $y$, respectively (notice that these paths could be degenerate, i.e. they could consist of a single vertex). Then
 $$
\r(x,y)=\om(\g)=\om(\g_x)+\om(\dl)+\om(\g_y)\ge\r(x,z)+\r(z,w)+\r(w,y)
\ge\r(x,y),
 $$
therefore, the inequalities $\om(\g_x)\ge\r(x,z)$, $\om(\dl)\ge\r(z,w)$, and
$\om(\g_y)\ge\r(w,y)$ from the filling definition must hold as equalities, in particular, $\r(z,w)=\om(\dl)$, and hence, the path $\dl$ is exact. The proof is complete.
 \end {proof}

A non-boundary path is said to be {\em exact}, if it is contained in some exact boundary path. Notice that, due to~\ref{ass:boud_paths}, the previous definition can be expanded correctly to the boundary paths also. Thus, an {\bf arbitrary path} in a filling is said to be {\em exact}, if it is contained in a boundary path whose weight is equal to the distance between its ending points. An edge forming an exact one-edge path is said to be  {\em exact\/} also. A vertex of a filling is said to be  {\em exact}, if it is contained in an exact path.

A boundary path is said to be {\em irreducible}, if it does not contain boundary vertices except the ending ones. Assertion~\ref{ass:boud_paths} implies the next result.

 \begin {cor}\label{cor:irreduceble_exact}
An edge of a filling is exact, if and only if it is contained in some irreducible exact path.
 \end {cor}

 \subsection {Minimal Parametric Fillings} \label{subsec:param}

 \begin {ass}\label{ass:min_fill_given_type}
Let $\cG$ be a minimal parametric filling of a finite pseudo-metric space  $(M,\r)$. Then each non-degenerate edge of this filling is exact.
 \end {ass}

 \begin {proof}
If some non-degenerate edge $e$ of a filling $\cG$ is not exact, then, due to finiteness of the number of boundary paths passing through the $e$, the weight of $e$ can be decreased preserving the validity of the inequalities from the filling definition. Therefore, the initial filling is not minimal, q.e.d.
 \end {proof}

 \begin {rk}
A degenerate edge of a minimal parametric filling need not be contained in an exact path, generally speaking. To construct an example, we choose $(M,\r)$ consisting of the vertices of some rectangle $p_1p_2p_3p_4$ in the Euclidean plane. By $x$ we denote the length of the diagonal of the rectangle. Let $G$ be the binary tree with the boundary $M$, such that each two edges coming to the diagonal vertices are adjacent, namely, if $e_i$ stands for an edge of the tree $G$ incident to $p_i$, then the edge $e_j$ is adjacent with the $e_{j+2}$, $j=1,\,2$. Further, if $\cG=(G,\om)$ is an arbitrary filling, then $\om(e_j)+\om(e_{j+2})\ge|p_jp_{j+2}|=x$, therefore, $\om(\cG)\ge 2x$. The latter implies that if we put $\om(e_i)=x/2$, and $\om(f)=0$, where $f$ is the unique interior edge of $G$, then we obtain a minimal parametric filling $\cG$.  It is clear that the weight of any boundary path $\g$ in $\cG$ passing through  $f$ is equal to $x$. On the other hand, such path $\g$ join the vertices lying on a side of the rectangle. Therefore, the distance between the ending vertices of $\g$ is less than $x$, and hence $\g$ is not exact. Thus,  $f$ is not an exact edge.
 \end {rk}

 \begin {cor}\label{cor:vertex_exact}
Each vertex of a minimal parametric filling is exact.
 \end {cor}

 \begin {proof}
Assume the contrary, i.e. let $v$ be a vertex of a minimal parametric filling $\cG$, which is not exact. The latter means that all the boundary paths going through $v$ are not exact. In particular, the latter implies that $\cG$ contains non-degenerate edges.

Due to Assertion~\ref{ass:exact_project}, all the paths in the base $\cG_F$ of the filling $\cG$, passing through $v'=\pi_F(v)$ are not exact also. Since $\cG$ contains non-degenerate edges, the base $\cG_F$ contains an edge incident to $v'$. All the boundary paths going through this edge are not exact and, since  $\cG_F$ is a minimal parametric filling, we obtain a contradiction with Assertion~\ref{ass:min_fill_given_type}. The proof is complete.
 \end {proof}

 \begin {cor}\label{cor:edge}
Each one-edge boundary path in a minimal parametric filling $\cG=(G,\om)$ of a pseudo-metric space $\cM=(M,\r)$ is exact.
 \end {cor}

 \begin {proof}
Let $\g$ be a such path, and $e=w_1w_2$ be its unique edge. Since $\cG$ is a filling of $\cM$, we have $\om(e)=\om(\g)\ge\r(w_1,w_2)$, therefore, if $\om(e)=0$, then $\r(w_1,w_2)=0$ also, and Corollary is valid.

Now let $e$ be a non-degenerate edge. Due to Assertion~\ref{ass:min_fill_given_type}, the edge $e$ is contained in an exact path and, hence, due to Assertion~\ref{ass:boud_paths}, $\om(\g)=\r(w_1,w_2)$, q.e.d.
 \end {proof}

 \begin {cor}\label{cor:edge2}
Each boundary path in a minimal parametric filling $\cG$ of a pseudo-metric space $\cM$, containing exactly one non-degenerate, edge is exact.
 \end {cor}

 \begin {proof}
Let $\g$ be a such path, and $F=\degen(\cG)$. Due to Assertion~\ref{ass:false_gen1}, the base $\cG_F$ of the graph $\cG$ is a minimal parametric filling of the space $\cM_F$. The path $\pi_F(\g)$ corresponding to $\g$ in $\cG_F$ is a boundary one-edge path. It remains to apply Corollary~\ref{cor:edge} and Assertion~\ref{ass:exact_project}. The proof is complete.
 \end {proof}

Let $v$ be an arbitrary vertex of a filling  $\cG$ of a metric space $\cM$. By $E_v$ we denote the set of all the edges from $\cG$, incident to $v$, and by $E_v^n$ we denote the subset of $E_v$, consisting of all non-degenerate edges. The family $N\ss E_v^n$ is said to be {\em non-exact}, if any pair of edges from $N$ is not contained in an exact path. By the{\em non-exactness index $\neind v$ of a vertex $v$} we call the maximal possible number of elements in a non-exact family $N\ss E_v^n$. Recall that the number of elements in the set $E_v$ is called by the {\em degree of the vertex $v$} and is denoted by  $\deg v$.

 \begin {ass}\label{ass:non_exact_vertex}
Let $v$ be an arbitrary interior vertex of a minimal parametric filling. Then $2\neind v\le\deg v$.  In the other words, each subset of $E_v$ consisting of more than $(\deg v)/2$ non-degenerate edges contains at least two edges belonging to an exact boundary path.
 \end {ass}

 \begin {proof}
Assume the contrary, i.e. there exists an interior vertex $v$ such that $\neind v>\deg v-\neind v$. Consider a non-exact family $N\ss E_v^n$ consisting of
$\neind v$ elements. For every pair $e_p$ and $e_q$ of edges from $N$, we consider all the boundary paths  $\g$ passing through $e_p$ and $e_q$, and by $\e_{pq}$ we denote the minimal difference between the weights of the paths $\g$ and the distances between their ending points considered as the points of the space $\cM$. Due to the definition of $N$, we have $\e_{pq}>0$. Now, choose an arbitrary positive $\e$ which is less than all the $\e_{pq}$ and all the weights of the edges from $N$. Change the weight function of the tree $\cG$, decreasing the weights of the edges from $N$ by $\e/2$ and increasing the edges from $E_v\sm N$ by $\e/2$. The tree remains a filling, but since $\neind v>\deg v-\neind v$, the total weight of the tree decreases. The latter contradicts to the minimality of $\cG$. The proof is complete.
 \end {proof}

An interior vertex $v$ of a filling is said to be {\em non-exact}, if $2\neind v>\deg v$.  Now Assertion~\ref{ass:non_exact_vertex} can be reformulated as follows: {\em minimal parametric fillings do not contain non-exact vertices}.

Let $G=(V,E)$ be an arbitrary tree. Let $v\in V$ be an interior vertex of degree $(k+1)\ge3$ adjacent with $k$ vertices $w_1,\ldots,w_k$ from $\d G$. Then the set of the vertices $\{w_1,\ldots,w_k\}$, and also the set of the edges  $\{vw_1,\ldots,vw_k\}$, are referred as {\em mustaches}. The number $k$ is called by the {\em degree}, and the vertex $v$ is called by the {\em common vertex of the mustaches}. An edge incident to $v$ and not belonging to $\{vw_1,\ldots,vw_k\}$ is called the {\em root edge\/} of the mustaches under consideration.

 \begin {ass}\label{ass:moustashes_param_gen}
Let $\cG=(G,\om)$ be a minimal parametric filling of a finite pseudo-metric space $\cM=(M,\r)$, and let $W=\{w_1,\ldots,w_k\}$ be a mustaches with common vertex $v$. Then some points $w_p,\,w_q\in W$ are joined by an exact path.
 \end {ass}

 \begin {proof}
Assume the contrary, i.e\. no two points from $W$ are joined by an exact path. Then all the edges  $vw_p$ are non-degenerate, since otherwise some vertices from $W$ would be connected either by a degenerate and, so, by an exact path, or by a path consisting of a single non-degenerate edge and, so, being exact in accordance with Corollary~\ref{cor:edge2}. But then $\neind v\ge k$, and therefore
 $$
2\neind v\ge 2k>k+1=\deg v,
 $$
where the strict inequality is valid since $k\ge 2$ in the definition of mustaches. Thus, the vertex $v$ is non-exact, that contradicts to assertion~\ref{ass:non_exact_vertex}. Proof is complete.
 \end {proof}

 \begin {cor}\label{cor:moustashes_param}
Let $\cG=(G,\om)$ be a minimal parametric filling of a pseudo-metric space $\cM=(M,\r)$. Then each its mustaches of degree $2$ are joined by an exact path.
 \end {cor}

 \begin {cor}\label{cor:three-joint}
let $v$ be an interior vertex of degree $m$ in a minimal parametric filling of some pseudo-metric space, which is incident to non-degenerate edges only. Then each set consisting of at most $m/2$ of these edges contains at least two edges belonging to an exact path. In particular, if the degree of $v$ is equal to $3$, then any pair of edges incident to $v$ belongs to an exact path.
 \end {cor}

The necxt result can be proved similarly to Assertion~\ref{ass:moustashes_param_gen}.

 \begin {cor}\label{cor:three-joint-one-degen}
Let $v$ be an interior vertex of degree $m$ of a minimal parametric filling of some pseudo-metric space, which is incident to  $(m-1)$ non-degenerate edges. Then there exists a pair of non-degenerate edges incident to  $v$, belonging to an exact path.
 \end {cor}

 \begin {cor}\label{cor:4-configuration}
Let a non-degenerate edge $e$ of a minimal parametric filling be incident to a pair of its vertices $v$ and $w$ of degree three. By $e_1$ and $e_2$ we denote the remaining edges incident to $v$, and by $f_1$ and $f_2$ we denote the remaining edges incident to  $w$.  Assume that all these four edges are also non0-degenerate. Then there exist two non-intersecting pairs of edges $(e_i,f_p)$ and $(e_j,f_q)$, $\{i,j\}=\{p,q\}=\{1,2\}$, such that each of them belong to an exact path.
 \end {cor}

 \begin {proof}
Indeed, in accordance with Corollary~\ref{cor:three-joint}, the pair $(e_1,e)$ belongs to an exact path. To be definite, assume that this path passes through the edge $f_1$ also. Du to the same Corollary~\ref{cor:three-joint}, there is an exact path passing through the edges  $(f_2,e)$. If this path contains $e_2$, then the pairs
$(e_1,f_1)$ and $(e_2,f_2)$ are the sought for. Otherwise, consider the path going through $(e_2,e)$. If it contains $f_2$, then the same pairs fit our conditions, and if it contains $f_1$, then we take the pairs $(e_2,f_1)$ and $(e_1,f_2)$. Corollary is proved.
 \end {proof}

 \subsection {Minimal fillings}
As we have already mentioned above, each minimal filling is a minimal parametric filling for some parametrization, therefore all the results of Subsection~\ref{subsec:param} remain valid for minimal fillings. In the present Subsection we list several more strong properties which hold for minimal fillings only.

Let $G$ be an arbitrary tree, and $e_1$ and $e_2$ be some its edges. By $\g$ we denote the unique path in $G$, which the edges $e_i$ are ending for. We say that an edge $f$ {\em lies between $e_1$ and $e_2$}, since it differs from the edges $e_i$ and belongs to the path $\g$. We also say, that a vertex  $v$ {\em lies between $e_1$ and $e_2$}, if $v$ is a non-ending vertex of the path $\g$.

 \begin {ass}\label{ass:boundary_nongen}
Let $\cG$ be a minimal filling of a pseudo-metric space $\cM$, and $e$ and $f$ be some edges from $\cG$. Assume that all the edges lying between  $e$ and $f$ are degenerate. Then the path joining $e$ and $f$ in $\cG$ is exact.
 \end {ass}

 \begin {proof}
Consider path $\bg$ containing $e$ and $f$, maximal with respect to inclusion. Due to our agreement concerning the types of filling under consideration, the path $\bg$ is a boundary path. If all the edges in $\bg$ are degenerate, then $\bg$ is exact and Assertion holds.

Now let $\bg$ contain non-degenerate edges. Orient $\bg$ such that $e$ precedes $f$. Assume firstly that all the edges of the path $\bg$ lying either after $e$, or before $f$ are degenerate. Without loss of generality, consider the first possibility. By $g$ we denote the last non-degenerate edge of the path $\bg$, and let $\cG_F$ be the vase of the filling $\cG$. Then $\pi_F(g)$ is a boundary edge in $\cG_F$. Assertion~\ref{ass:min_fill_given_type} and Corollary~\ref{cor:irreduceble_exact} imply that $\pi_F(g)$ belongs to some irreducible exact path $\g'$. Let $w$ be the end vertex of the path $\bg$. Notice that $\pi_F(w)$ is a boundary vertex of the filling
$\cG_F$, incident to the edge $\pi_F(g)$. Since $\g'$ is an irreducible path, then  $\pi_F(w)$ is its ending vertex. Let $u'$ be the ending vertex of the path $\g'$, distinct from $\pi_F(w)$, and $u$ be an arbitrary boundary vertex from $\cG$, such that $\pi_F(u)=u'$. Clearly, that the path $\g$ joining $u$ and $w$  in $\cG$ contains $e$ and $f$, and moreover, $\g$ is exact due to exactness of $\g'$. Thus, in this case Assertion is proved.

At last, assume that the path $\bg$ contains non-degenerate edges as after $e$, so as before $f$. By $e_1$ we denote the last non-degenerate edge on $\bg$, located before $f$, and by $e_2$ we denote the first non-degenerate edge lying after $e$. Assume that in this case Assertion is not valid. By $v$ we denote an arbitrary vertex lying between $e_1$ and $e_2$. Put $e_i'=\pi_F(e_i)$, $v'=\pi_F(v)$.  Due to Assertion~\ref{ass:exact_project}, each boundary path $\g'$ in $\cG_F$ passing through the both $e_i'$ is not exact. By $\e$ we denote an arbitrary positive number which is as less than the both weights $\om_F(e_i')=\om(e_i)$, so as less than each difference between the weight of such path  $\g'$ and the distance between its endpoints (die to our assumption, such $\e$ does exist).

Split the vertex $v'$ into an interior vertex $w_1'$ and a vertex $w_2'$ in such a way that $w_1'$ is incident exactly to $e_1'$, $e_2'$, and $w_1'w_2'$. Define the weight function $\om_F'$ on the resulting trees as follows. Put $\om_F'(w_1'w_2')=\e/2$, $\om_F'(e_i')=\om_F(e_i')-\e/2$,
and let $\om_F'=\om_F$ on all the remaining edges. It is clear that we obtain a filling of the space  $\cM_F$ again, but the weight of this filling is less, then the weight of the initial one. The latter contradicts to the minimality of $\cG_F$. The proof is complete.
 \end {proof}

 \begin {cor}\label{cor:any_pair}
Each path in a minimal filling, consisting of at most two edges is exact. In particular, each pair of adjacent edges of a minimal filling belongs to an exact boundary path.
 \end {cor}

A finite pseudo-metric space $(M,\r)$ is said to be {\em additive}, if there exists a weighted tree $(G,\om)$ such that $\d G=M\ss V(G)$ and $\r$ coincides with the pseudo-metric $d_\om$ generated by the weight function $\om$. Each such a weighted tree $(G,\om)$ is referred as  {\em generating\/} for the space $(M,\r)$.

Corollary~\ref{cor:any_pair} implies the next result.

 \begin {cor}\label{cor:star}
Assume that each boundary path of a minimal filling $\cG=(G,\om)$ of a finite pseudo-metric space $\cM=(M,\r)$ contains at most two non-degenerate edges, i.e\. the base $\cG_F$ of this filling has at most one vertex of degree more than $1$. Then $\cM$ is additive, and  $\cG$ is a generating tree for $\cM$.
 \end {cor}

 \begin {proof}
Let $p$ and $q$ be arbitrary points from $\cM$, such that $\r(p,q)\ne0$, and $\g$ be the path in $\cG$, joining $p$ and $q$. The path $\pi_F(\g)$ is a boundary one and by the assumptions it consists of at most two edges, hence, due to Corollary~\ref{cor:any_pair}, the path $\pi_F(\g)$ is exact. Therefore, $\g$ is exact and, hence, $d_\om(p,q)=\om(\g)=\r(p,q)$, q.e.d.
 \end {proof}

Let us enforce the previous Corollary. Let $\cG=(G,\om)$ be a weighted graph with a boundary $\d G$ and vertex set $V$. Choose an arbitrary vertex $v\in V$. If $\cG$ is non-degenerate, then by $\d_v$ we denote the set consisting of all the vertices from $\d G$, adjacent to $v$, together with the vertex  $v$, providing $v\in\d G$. In the general case, consider the base $\cG_F=(G_F,\om_F)$ of the graph $\cG$ and let $\pi_F\:V\to V_F$ be the natural projection, and $w=\pi_F(v)$. Then we put $\d_v=\pi_F^{-1}(\d_w)$. The next Assertion can be proved similarly to Corollary~\ref{cor:star}.

 \begin {cor}\label{cor:moustaches-additive-subspace}
For any vertex $v$ of a minimal filling $\cG$ of a pseudo-metric space $\cM=(M,\r)$, the subspace $(\d_v,\r)$ is additive. In particular, any mustaches of a minimal filling $\cG$ are an additive subspace in $\cM$.
 \end {cor}

In the next Assertion the exactness of several elements of a minimal filling is investigated.

 \begin {cor}\label{cor:any_edge_is_exact}
Let $\cG$ be a minimal filling consisting of at most two edges,  $v$ be its arbitrary vertex, and  $e$ be its arbitrary edge, degenerate or non-degenerate. Then
 \begin {description}
  \item[(1)] the vertex $v$ and the edge $e$ are exact\/\rom;
  \item[(2)] if $s$ stands for the degree of the vertex $v$ and $s\ge2$, then $v$ belongs to at least $C_s^2$ of exact paths, where $C_s^2$ is the binomial coefficient\/\rom;
  \item[(3)] if $m$ is the maximal degree of the vertices incident to $e$, then $e$ belongs to at least $(m-1)$ exact paths\/\rom;
  \item[(4)] moreover, as $v$, so as $e$ belongs to at least two exact paths.
 \end {description}
 \end {cor}

 \begin {proof}
(1) The exactness of the vertex $v$ follows from Corollary~\ref{cor:vertex_exact}, and the exactness of the edge $e$ --- from Corollary~\ref{cor:any_pair}.

(2) Due to Corollary~\ref{cor:any_pair}, each pair of edges incident to $v$, belongs to a boundary path. Distinct pairs belong to distinct paths.

(3) Let $e=vw$, and $m$ be the degree of the vertex $v$. Then $v$ is incident to $e$ and some edges $e_1,\ldots,e_{m-1}$ and, due to Corollary~\ref{cor:any_pair}, each pair $(e,e_i)$ belongs to an exact boundary path, and these paths are pairwise distinct.

(4) We start with the case of the vertex $v$. If its degree $s$ is more than $2$, then $C_s^2>2$ and Item~(4) follows from Item~(2). If $s=1$ or $2$, then, du to our agreement, $v$ is a boundary vertex.

At first, let the degree of $v$ be $2$, and $e_1$ and $e_2$ be the edges incident to $v$. In accordance with Item~(1), the both  $e_i$ are exact. Due to Corollary~\ref{cor:irreduceble_exact}, each $e_i$ belongs to an irreducible exact boundary path. But since $v$ is a boundary vertex, then these paths are distinct, and so,  $v$ belongs to at least two distinct exact boundary paths.

Now, let the degree of $v$ be equal to $1$. By $w$ we denote the unique vertex adjacent to $v$. If the degree of $w$ is more than two, then the edge $vw$
is adjacent with at least two distinct edges and Corollary~\ref{cor:any_pair} implies the required. If the degree of $w$ equals $2$, then $w$ is a boundary
vertex, and the edge $vw$ is an exact boundary path itself. The second boundary path can be found by applying Corollary~\ref{cor:any_pair} to the edge $vw$ and the second
edge incident to $w$.%
 \footnote{%
The first proof of this non-trivial statement was obtained by I.~Laut, a student of mechanical and mathematical faculty of Moscow State University from similar considerations.
 }

Now consider the edge $e$. Notice that since $\cG$ contains at least two edges, then $m\ge 2$. If $m\ge3$, then the required is proved. Now let
$m=2$. By $w$ we denote the degree $2$ vertex incident to $e$, and by $u$ we denote the remaining vertex of the edge $e$. By definition of $m$, the degree of
$u$ does not exceed $2$, so, due to our agrement, the both  vertices $u$ and $w$ are boundary ones. As above, the edge $e$ is an exact boundary path itself. The second exact boundary path can be found by applying Corollary~\ref{cor:any_pair} to the edge $e$ and the edge adjacent to $e$ by the vertex $w$.

The proof is complete.
 \end {proof}

Corollary~\ref{cor:any_pair} gives an opportunity to enforce Corollary~\ref{cor:4-configuration} as follows.

 \begin {cor}\label{cor:4-configuration-mf}
Let $e=vw$ be an interior edge of a minimal filling $\cG$ of a finite pseudo-metric space $\cM$ \(in accordance with our agreement, the degrees of the vertices $v$ and $w$ are more than $2$\). Let  $e_1$ and $e_2$ be some edges, which are distinct from $e$ and incident to $v$, and $f_1$ and $f_2$ be some edges, which are distinct from $e$ and incident to $w$. Then there exist two non intersecting pairs of edges $(e_i,f_p)$ and $(e_j,f_q)$, $\{i,j\}=\{p,q\}=\{1,2\}$, each of which belongs to an exact path.
 \end {cor}

 \begin {proof}
Indeed, due to Corollary~\ref{cor:any_pair}, the pair $(e_1,e)$ belongs to an exact path. To be definite, assume that this path passes through the edge  $f_1$ also. Due to the same Corollary~\ref{cor:any_pair}, the edges $(f_2,e)$ belong to an exact path also. If this latter path passes trough $e_2$, then the pairs $(e_1,f_1)$ and $(e_2,f_2)$ are required. Otherwise, we consider an exact path containing $(e_2,e)$. If it passes through $f_2$, then we can tale the same pairs of edges, and if it passes trough $f_1$, then the required pairs are $(e_2,f_1)$ and $(e_1,f_2)$. Corollary is proved.
 \end {proof}

Recall that a triangle $uvw$ in a pseudo-metric space is said to be {\em non-degenerate}, if all its three triangle inequalities hold strictly. Otherwise the triangle is said to be  {\em degenerate}.  a point $p$ of a finite pseudo-metric space $(M,\r)$ is called {\em degenerate}, if there exist points $q$ and $r$ in $M$, which are distinct from $p$ and such that $\r(q,r)=\r(q,p)+\r(p,r)$.

 \begin {cor}
Let $\cG$ be a minimal filling of a finite pseudo-metric space $\cM$, and $e$ and $f$ be some edges from $\cG$. assume that all the edges lying between $e$ and $f$ are degenerate, and that some boundary vertex $v$ lies between $e$ and $f$. Then  $v$ is a degenerate point of the space $\cM$.
 \end {cor}

 \begin {proof}
In accordance with Assertion~\ref{ass:boundary_nongen}, there exists an exact boundary path $\g$ joining $e$ and $f$. By $v_1$ and $v_2$ we denote the ending vertices of the path $\g$, and let $\g_i$ be the part of the path $\g$ between $v_i$ и $v$. Due to Corollary~\ref{cor:edge2}, the paths $\g_i$ are also exact. Therefore,
 $$
\r(v_1,v_2)=\om(\g)=\om(\g_1)+\om(\g_2)=\r(v_1,v)+\r(v,v_2),
 $$
q.e.d.
 \end {proof}

A pseudo-metric space is said to be {\em non-degenerate}, if all triangle inequalities in it are strict. Notice that a non-degenerate pseudo-metric space is a metric space.

 \begin {cor}\label{cor:boundary-vertices-nondegenerate}
Let $\cG$ be a minimal filling of a finite pseudo-metric space $\cM$. Then each boundary vertex of degree more than $1$ is a degenerate point of the space $\cM$. In particular, if the space $\cM$ is non-degenerate, then the degree of any boundary vertex of the tree $\cG$ is equal to $1$ and all the boundary edges are non-degenerate.
 \end {cor}

 \begin {cor}
A minimal filling $\cG$ of a finite non-degenerate non-additive metric space  $\cM$ contains a non-degenerate interior edge.
 \end {cor}

 \begin {proof}
Consider the base $\cG_F$ of the filling $\cG$. Due to Corollary~\ref{cor:boundary-vertices-nondegenerate}, all the boundary vertices of the filling  $\cG_F$ have degree $1$, therefore, if $\cG$ does not have non-degenerate interior edges, then the graph $\cG_F$ has unique vertex of degree more than $1$, and, due to Corollary~\ref{cor:star}, the space $\cM$ is additive, a contradiction. The proof is complete.
 \end {proof}

Notice that the set of all pseudo-metric spaces consisting of $n$ points can be naturally identified with a convex cone in $\R^{n(n-1)/2}$ (it suffices to enumerate the set of all two-elements subsets of these spaces and assign to each such space the vector of the distances between the points of the pairs). This representation  gives us an opportunity to speak about topological properties of families of metric spaces consisting of a fixed number of points.

We say, that some property holds for a {\em pseudo-metric space in general position}, if for any $n$ this property is valid for an everywhere dense set of  $n$-point pseudo-metric spaces.

 \begin {conj}
A minimal filling of a finite metric space in general position is a non-degenerate binary tree.
 \end {conj}

Let $\cG$ be a minimal filling of a finite pseudo-metric space $\cM=(M,\r)$. Define on $M$ the {\em graph of exact paths for $\cG$}, joining points $p$ and $q$ from $M$ by an edge, if and only if they are joined by an exact boundary path in $\cG$.

 \begin {ass}
The graph of exact paths for a minimal filling $\cG$ is connected, i.e\. for any partition $\{M_1, M_2\}$ of the set $M$ there exist points $p\in M_1$ and $q\in M_2$, joined by an exact path in $\cG$.
 \end {ass}

 \begin {proof}
Assume the contrary, i.e\. let there exist a partition  $\{M_1, M_2\}$ of the set $M$, such that for any two points $p\in M_1$ and $q\in M_2$, the path joining them in $\cG$ is not exact.

Consider the subgraph $G_i$ of the graph $G$ generated by the edges of all the exact boundary paths joining the vertices from the set $M_i$. Due to the assumption, each exact boundary path in $\cG$ joins the vertices from the same $M_i$, therefore, due to Corollary~\ref{cor:any_edge_is_exact}, we have: $G=G_1\cup G_2$. Since $\cG$ is a connected graph (a tree), the subgraphs $G_i$ have a non-empty intersection which is denoted by $H$.

Let $v$ be a vertex of $H$, whose degree in the forest $H$ does not exceed $1$. By $\g_i$ we denote some exact boundary path joining the vertices from  $M_i$ and containing $v$. Notice that $v$ is an interior vertex of the path  $\g_i$, therefore its degree in $\g_i$ is equal to $2$. Since the degree of the vertex $v$ in $H$ does not exceed $1$, then $\g_i$ contains an edge $e_i$ which is incident to $v$ and is not contained in the graph $H$, in particular, $e_1\ne e_2$. Moreover, $e_1$ and $e_2$ can not belong to an exact boundary path, since otherwise at least one of these edges belongs to $H$. But the latter contradicts to Corollary~\ref{cor:any_pair}. The proof is complete.
 \end {proof}

 \section {Tours, Perimeters, and Minimal Fillings} \label{sec:tours}
 \markright {\thesection.~Tours, Perimeters, and Minimal Fillings.}
Let  $G=(V,E)$  be an arbitrary tree with a boundary  $M$. Recall that in accordance with our Agreement, the set $M$ contains all the degree $1$ and $2$ vertices of the tree  $G$.  Remove some edge $e$ from $G$, and let $G_1$ and $G_2$ be the connected components of the resulting forest. We put $M_i=M\cap G_i$.

 \begin {lem}
Under the above notations, the sets $M_i$ are not empty.
 \end {lem}

 \begin {proof}
Indeed, if  $G_i$ consists of a single vertex $v$, then the degree of $v$  in $G$ is equal to $1$, and, hence, $v\in M_i$. If the tree $G_i$ contains an edge, then
$G_i$ contains at least two vertices of degree $1$, and only one of them is incident to the edge  $e$ in the tree $G$, therefore the remaining ones have degree  $1$ in $G$ and, hence, belong to $M_i$. Lemma is proved.
 \end {proof}

By $\cP_G(e)$ we denote the partition $\{M_1,M_2\}$ of the set $M$ obtained.

Let $S$ be a finite set. By a  {\em cyclic order on the set  $S$} we call an arbitrary cyclic permutation $\pi\:S\to S$. Two elements from  $S$ are said to be  {\em adjacent\/} with respect to a cyclic order $\pi$, if one of them is the  $\pi$-image of the other one. An enumeration $(s_1,\ldots,s_k)$ of the elements from  $S$ is said to be {\em compatible\/} with a cyclic order $\pi$, if $\pi(s_i)=s_{i+1}$ for any $i$, $i<k$. It is clear, that an enumeration  $(s_1,\ldots,s_k)$, $k=|S|$, is compatible with a cyclic order $\pi$, if and only if $s_{i+1}=\pi^i(s_1)$ for all $i$, $i<k$. For any cyclic order on the set $S$ there exist exactly $k$ compatible enumerations.

A cyclic order $\pi$ on $M$ is said to be {\em planar with respect to $G$} or is called a {\em tour of $G$\/}, for any $e\in E$ and for each $M_i\in\cP_G(e)$ there exists unique point $p\in M_i$, such that $\pi(p)\not\in M_i$. The latter means that there exists an enumeration of the set $M$ compatible with $\pi$, such that
the elements from the set $M_1$ precede the elements from the set $M_2$ in it.

We give an equivalent definition of the planar order in terms of embeddings. Let $G'$ be some embedding of the tree  $G$ into the plane. Consider a walk around the tree $G'\ss\R^2$. We draw the points of $M$ consecutive with respect to this walk as a consecutive points of the circle $S^1$. Notice that each vertex $p$ from $M$ appears $\deg p$ times.  For each vertex $p\in M$ of degree more than $1$, we choose one arbitrary point from the corresponding points of the circle. So, we construct an injection $\nu\:M\to S^1$. Define cyclic permutation $\pi$ as follows:  $\pi(p)=q$, where $\nu(q)$ follows after $\nu(p)$ on the circle $S^1$. We say that the cyclic order $\pi$ {\em is generated by the embedding $G'$}. It is clear that the embedding  $G'$ generates  $2\prod_{p\in M}\deg p$ of the cyclic orders.

 \begin {ass}
A cyclic order on $M$ generated by an embedding $G'$ of the tree  $G$, is planar with respect to $G$. Inversely, each planar order on  $M$ with respect to  $G$ is generated by some embedding  $G'$ of the tree  $G$.
 \end {ass}

 \begin {proof}
The first statement is evident since, for any edge $e$ of the tree $G$, the sets $\nu(M_i)$, $M_i\in\cP_G(e)$, do not alternate with each other on the circle $S^1$.

Let us prove the inverse statement using induction over the number $m$ of edges of the tree  $G$. If $m=1$, then $M$ consists of two points and the unique cyclic order on  $M$ is generated by the unique embedding of the single-edge tree.

Now, let $m>1$. Since $M$ contains all the vertices of degree  $1$ and $2$, it consists of at least three points. Consider an arbitrary order $\pi$ on $M$, which is planar with respect to $G$, and let $x$ be an arbitrary degree $1$ vertex of $G$. By $\ti G$ we denote the tree obtained from the tree $G$ by removing the vertex  $x$ and the unique edge $e=xy$ incident to $x$.

By $\ti M$ we denote the set $M\sm\{x\}$.  It is clear that the set $\ti M$ consists of the vertices of the tree $\ti G$ and contains all its vertices of degree $1$. Consider the mapping $\ti\pi$ defined on $\ti M$ as follows: $\ti\pi\bigl(\pi^{-1}(x)\bigr)=\pi(x)$ and $\ti\pi=\pi$ for the remaining points from $\ti M$. Since $M$ contains at least three points, the mapping $\ti\pi$ is a well-defined cyclic order. Let us show that the order $\ti\pi$ is planar with respect to $\ti G$.

Notice that each edge  $\ti e$ of the tree $\ti G$ is an edge of the tree $G$ also. Let $\{M_1,M_2\}=\cP_G(\ti e)$, and let $(p_1,\ldots,p_n)$ be an enumeration of the elements from $M$, which is compatible with the planar order $\pi$ an such that  $M_1=\{p_1,\ldots, p_k\}$, and $M_2=\{p_{k+1},\ldots,p_n\}$. Without loss of generality assume that  $x\in M_2$, i.e\. $x=p_i$, $i\ge k+1$. Since $\ti e\ne xy$, then $k\le n-2$ in this case. Consider the enumeration of the set $\ti M$, which is obtained by removing of $x$, namely, $(p_1,\ldots,p_{i-1},p_{i+1},\ldots,p_n)$. Then $\ti\pi(p_j)=p_{j+1}$ for all $j\ne i-1$, and $\ti\pi(p_{i-1})=\pi(x)$, where $\pi(x)=p_{i+1}$ if $i<n$, and $\pi(x)=p_1$ if $i=n$. At last, if  $\{\ti M_1,\ti M_2\}=\cP_{\ti G}(\ti e)$, then $\ti M_1=M_1$ and $\ti M_2=M_2\sm\{x\}$. Therefore, the enumeration of the set $\ti M$ constructed above is compatible with the cyclic order $\ti\pi$, and the elements of the set $\ti M_1$ precede the elements of the set  $\ti M_2$ in this enumeration. Due to arbitrariness of the edge $\ti e$, the latter implies that the cyclic order  $\ti\pi$ is planar with respect to the tree $\ti G$.

Due to the inductive assumptions, the planar order $\ti\pi$ is generated by some embedding  $\ti G'$ of the tree $\ti G$. It remains to reconstruct the embedding  $\ti G'$ to an embedding  $G'$ of the tree $G$ by emitting some segment from the image of the vertex $y$. This segment corresponds to the edge  $e=xy$. We proceed as follows.

If $y\in\ti M$, then we emit the segment on the same side of the tree $\ti G$ as we walk around it from $\ti\pi^{-1}(y)$ to $\ti\pi(y)$. Otherwise, we emit the segment on the same side of the tree $\ti G'$ as we walk around it from $\pi^{-1}(x)$ to $\pi(x)$. It remains to notice that the order generated on $M$ by the embedding $G'$ of the tree $G$ constructed just now coincides with $\pi$.  Assertion is proved.
 \end {proof}

 \begin {ass}\label{ass:2weights}
Let $\cG=(G,\om)$ be a weighted tree with a boundary $M$ and $\pi$ be an arbitrary cyclic order on $M$. Then
 $$
\sum_{p\in M}d_\om\bigl(p,\pi(p)\bigr)\ge2\om(G).
 $$
Moreover, the equality is achieved, if and only if $\pi$ is a tour around $G$.
 \end {ass}

 \begin {proof}
By $\g_p$ we denote the path in the tree $G$ joining $p$ and $\pi(p)$. Then $d_\om\bigl(p,\pi(p)\bigr)=\om(\g_p)$. To prove the inequality it suffices to verify that for any edge  $e$ from $G$ there exist distinct points $p$ and $q$ from $M$, such that the paths $\g_p$ and $\g_q$ contain $e$. To do that, consider the partition $\cP_G(e)=\{M_1,M_2\}$. It is easy to see, hat each $M_i$ contains a point  $p_i$ such that  $\pi(p_i)$ belongs to $M_j$, $j\ne i$. We can take $p_1$ and $p_2$ as the $p$ and $q$.  It remains to notice that the cyclic order $\pi$ is planar, if and only if each edge  $e$ is contained in exactly two paths $\g_p$. The proof is complete.
 \end {proof}

Let $\cM=(M,\r)$ be a finite pseudo-metric space, and $\pi$ be an arbitrary cyclic order on $M$.  The {\em perimeter of the space $\cM$ with respect to the order  $\pi$} is the value
 $$
P(\cM,\pi)=\sum_{p\in M}\r\bigl(p,\pi(p)\bigr),
 $$
and  $\min_\pi P(\cM,\pi)$, where minimum is taken over all possible cyclic orders $\pi$ on $M$ is called the  {\em perimeter of the pseudo-metric space $M$} and is denoted by $P(\cM)$. Notice that $P(\cM)$ is the length of the path solving the traveling salesman problem, see~\cite{PrepSham}.

Further, the {\em half-perimeter $p(\cM,\pi)$ of the space $\cM$ with respect to the order $\pi$} and simply the {\em half-perimeter $p(\cM)$} of a finite pseudo-metric space  $\cM$ is defined as the half of the corresponding perimeter:
 $$
p(\cM,\pi)=P(\cM,\pi)/2, \qquad p(\cM)=P(\cM)/2.
 $$

Besides, in what follows we need the following notations. If $\cG=(G,\om)$ is some filling of a finite pseudo-metric space $\cM=(M,\r)$, then the set of all cyclic orders on $M$, which are planar with respect to a tree $G$, i.e. the set of all tours of  $G$, is denoted by  $\cO(G)$ or by  $\cO(\cG)$. We also say that each such tour (cyclic order) is defined on  $\cM$.

Assertion~\ref{ass:2weights} and definition of the filling imply the following result.

 \begin {cor}\label{cor:fill_ge_halfp}
Let $\cG=(G,\om)$ be an arbitrary filling of a finite pseudo-metric space $\cM=(M,\r)$ and $\pi\in\cO(G)$. Then $\om(G)\ge p(\cM,\pi)$, in other words,
 $$
\om(G)\ge\max_{\pi\in\cO(G)}p(\cM,\pi)
 $$
and, hence,
 $$
\mf(\cM)\ge\min_G\max_{\pi\in\cO(G)}p(\cM,\pi)\ge p(\cM),
 $$
where minimum can taken over as over all the trees with the boundary $M$, so as over all the binary trees with the boundary $M$ only.

Further, if
 $$
\om(G)=\max_{\pi\in\cO(G)}p(\cM,\pi),
 $$
then $\cG$ is a minimal parametric filling, and
 $$
\mpf(\cM,G)=\max_{\pi\in\cO(G)}p(\cM,\pi).
 $$
 \end {cor}

Let $\cG=(G,\om)$ be a filling of a finite pseudo-metric space  $\cM=(M,\r)$. A tour $\pi\in\cO(G)$ is said to be {\em exact}, if for any point $p\in M$ the equality
$\r\bigl(p,\pi(p)\bigr)=d_\om\bigl(p,\pi(p)\bigr)$ holds, in other words, if each boundary path joining boundary points consecutive with respect to $\pi$ is exact.

 \begin {ass}\label{ass:exact_order}
Let $\cG=(G,\om)$ be a filling of a finite pseudo-metric space  $\cM=(M,\r)$, and  $\pi\in\cO(G)$ be a tour of $G$.  Then the tour $\pi$ is exact, if and only if  $\om(G)=p(\cM,\pi)$. Thus, if $\pi$ is an exact tour of $G$, then $\cG$ is a minimal parametric filling, and $\mpf(\cM,G)=p(\cM,\pi)=\max_{\s\in\cO(G)}p(\cM,\s)$.
 \end {ass}

 \begin {proof}
Let $\pi$ be an exact tour, then, due to Assertion~\ref{ass:2weights}, we have
 $$
2\om(G)=
\sum_{p\in M}d_\om\bigl(p,\pi(p)\bigr)=
\sum_{p\in M}\r\bigl(p,\pi(p)\bigr)=
2p(\cM,\pi).
 $$
Inversely, let $\om(G)=p(\cM,\pi)$. Due to the filling and tour definitions, we have:
 $$
2\om(G)=
\sum_{p\in M}d_\om\bigl(p,\pi(p)\bigr)\ge
\sum_{p\in M}\r\bigl(p,\pi(p)\bigr)=
2p(\cM,\pi)=2\om(G),
 $$
and hence, $\r\bigl(p,\pi(p)\bigr)=d_\om\bigl(p,\pi(p)\bigr)$ for all $p\in M$, and so, the tour $\pi$ is exact.

The remaining part of Assertion follows from Corollary~\ref{cor:fill_ge_halfp}. The proof is complete.
 \end {proof}

Thus, Assertion~\ref{ass:exact_order} gives an opportunity to calculate the weight of a minimal parametric filling as the maximum of the half-perimeters of all the tours of this filling, but under the assumption that an exact tour exists. If there is no an exact tour, then combining Assertion~\ref{ass:exact_order} and Corollary~\ref{cor:fill_ge_halfp} we conclude that $\mpf(\cM,G)>\max_{\s\in\cO(G)}p(\cM,\s)$.

Let us find out, if a minimal parametric filling always possesses an exact tour. As it is mentioned in the proof of Assertion~\ref{ass:2weights}, for each edge $e$ of the tree  $G$ and for any cyclic order $\pi$ on $\d G$, there exists a path going through  $e$ and joining boundary vertices consecutive with respect to $\pi$. This observation leads to the following result.

 \begin {cor}\label{cor:non-exact-edge-round-trip}
If a filling contains non-exact edge, then no one tour is exact.
 \end {cor}

Recall that above, in Remark from Subsection~\ref{subsec:param}, we give the example of a minimal parametric filling $\cG$ of a four-point metric space and find out that $\cG$ contains non-exact edge of zero weight. So, Corollary~\ref{cor:non-exact-edge-round-trip} implies that this filling does not possess an exact tour. Thus, the formula $\mpf(\cM,G)=\max_{\s\in\cO(G)}p(\cM,\s)$ is not valid in general case.

 \begin {prb}
Find out, which fillings possess an exact tour.
 \end {prb}

We have shown above that the presence of a non-exact edge is an obstacle to an exact tour existence. But there exist minimal parametric fillings all whose edges are exact, but without an exact tour.

 \begin {examp}
This example is obtained by Z.~Ovsjannikov, a student of mechanical and mathematical department of Moscow state university, by means of {\em Mathematica\/} software. Consider the metric space $\cM=(M,\r)$, where $M=\{p_1,p_2,p_3,p_4,p_5\}$ and the distance function $\r$ is given by the following matrix:
 $$
 \left(
\begin {array}{ccccc}
 0&4&6&6&4 \\
 4&0&4&6&6 \\
 6&4&0&4&6 \\
 6&6&4&0&4 \\
 4&6&6&4&0
 \end {array}\right).
 $$
This space is conventionally shown in Figure~\ref{fig:examp5} as a complete graph, and solid lines joins the pairs of points such that the distance between them is  $6$, and dashed lines correspond to distance  $4$. As a topology of minimal parametric filling  $\cG=(G,\om)$ we choose the binary tree with mustaches $\{p_1,p_2\}$ and
$\{p_3,p_5\}$.

\Pic{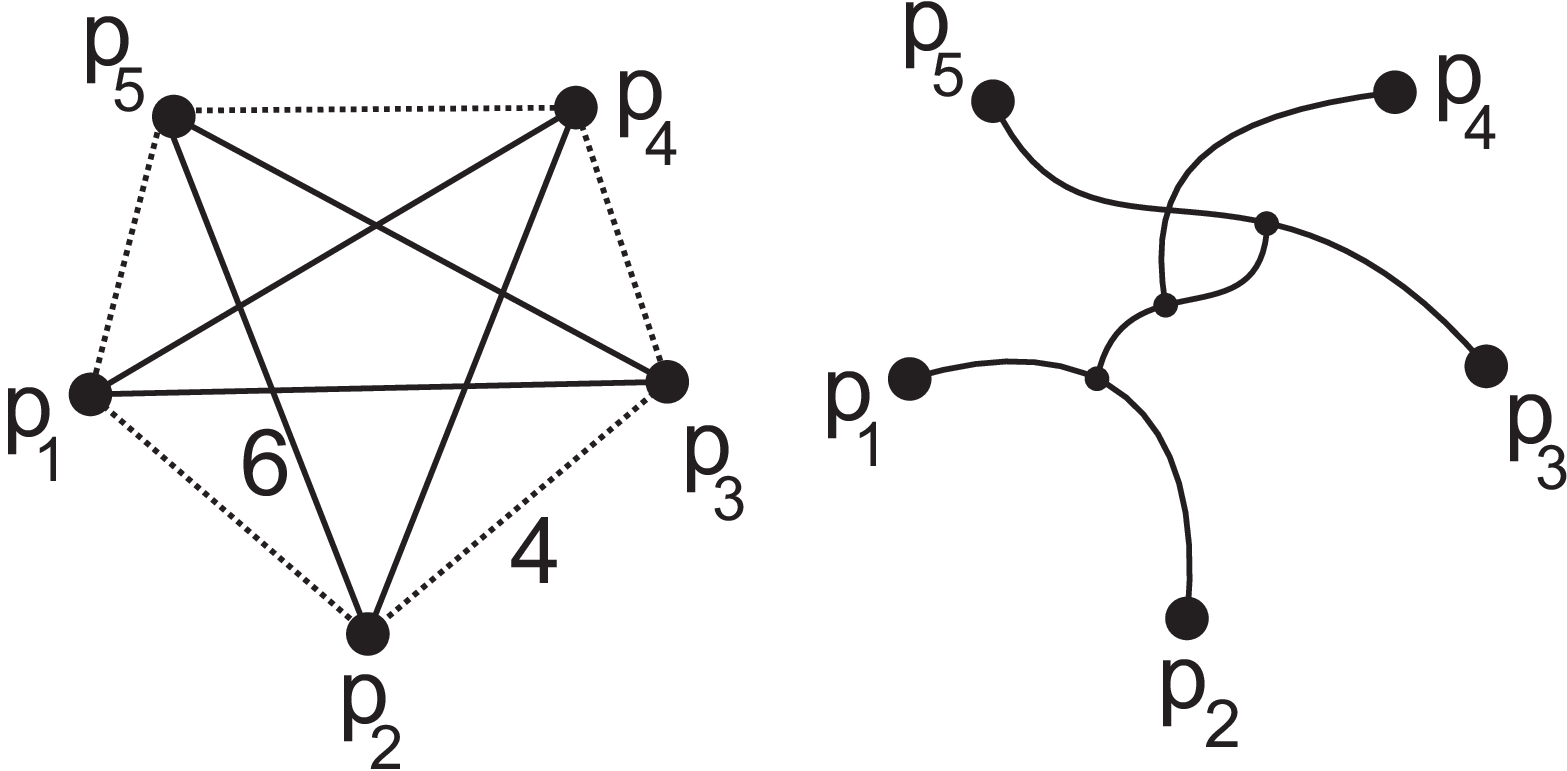}{This topology on this metric space generates minimal parametric filling all whose edges are exact, but which does not have an exact tour. Therefore, its maximal half-perimeter is less than the weight of the filling.}{fig:examp5}{200}

By  $x_i$ we denote the weight of the boundary edge of the tree $\cG$, incident to $p_i$, and by  $y_1$ and $y_2$ we denote the weights of the interior edges, the first one is a root-edge for the mustaches $\{p_1,p_2\}$, and the second one is a root-edge for the mustaches $\{p_3,p_5\}$.  Then the definition of filling and Corollary~\ref{cor:moustashes_param} imply that the weights of a minimal filling meet the following relations:
 \begin {gather*}
x_1+x_2=4,\qquad x_3+x_5=6, \\
x_1+y_1+x_4\ge6, \qquad x_2+y_1+x_4\ge6,
 \end {gather*}
and so, excluding $x_2=4-x_1$ from the second inequality, we get
 $$
x_1+y_1+x_4\ge6, \qquad -x_1+y_1+x_4\ge2.
 $$
Summing the later inequalities, we conclude that $y_1+x_4\ge4$.  On the other hand, the weight of the filling  $G$ is equal to $\sum_ix_i+y_1+y_2=10+(x_4+y_1)+y_2\ge14+y_2\ge14$. At last, it is not difficult to verify that the weight function $\om$ defined as follows: $x_5=4$, $y_2=0$, and
$x_1=x_2=x_3=x_4=y_1=2$ gives a filling $\cG=(G,\om)$ of the space $M$, and $\om(G)=14$.  Therefore,  $\cG$ is a minimal parametric filling.

Notice that all the edges from $\cG$ are exact. Indeed, the edges incident to mustaches are exact automatically. Further, the weight of the path $\g_{34}$ joining $p_3$ and$p_4$ is equal to $4$, hence $\g_{34}$ is exact. So, the edge incident to  $p_4$ and the root-edge of the mustaches $\{p_3,p_5\}$ are exact. At last, the weight of the path $\g_{24}$ joining $p_2$ and $p_4$ is equal to $6$, hence $\g_{24}$ is exact and the root-edge of the mustaches $\{p_1,p_2\}$ is also exact.

The set $\cO(G)$ consists of the next four planar orders (up to the inverse orders):
 $$
\pi_1=(1,2,3,5,4),\quad\pi_2=(1,2,5,3,4),\quad\pi_3=(2,1,3,5,4),
\quad\pi_4=(2,1,5,3,4).
 $$
The corresponding half-perimeters are as follows:
 $$
p(\cM,\pi_1)=12,\quad p(\cM,\pi_2)=13,\quad p(\cM,\pi_3)=13,\quad
p(\cM,\pi_4)=12,
 $$
therefore, the maximal half-perimeter in this example is less than the weight of minimal parametric filling, and there is no an exact tour.
 \end {examp}

Notice that all the above examples of minimal parametric fillings without exact tours contain degenerate edges. It turns out that it is also possible to give such example of non-degenerate filling.

 \begin {examp}
In the previous example, let us pass to the base $\cG_F$, where $F$ consists of the unique degenerate edge of the filling  $\cG$. Since this edge is interior, the base $\cM_F$ of the space $\cM$ coincides with $\cM$. As a result, we obtain a non-degenerate minimal parametric filling $\cG_F=(G_F,\om_F)$ of the same space $\cM$ of the same weight $14$.  All the edges of $\cG_F$ are exact, since the weights of all the paths are preserved under the projection onto the base. The number of tours increases, since the three boundary edges incident to the degree $4$ vertex cab be ordered in an arbitrary way. Namely, the set $\cO(G_F)$ consists of the following $6$ tours (up to the inverse orders):
 \begin {gather*}
\pi_1=(1,2,3,4,5),\qquad\pi_2=(1,2,3,5,4),\qquad\pi_3=(1,2,4,3,5)\\
\pi_4=(1,2,4,5,3),\qquad\pi_5=(1,2,5,3,4),\qquad\pi_6=(1,2,5,4,3).
 \end {gather*}
The corresponding half-perimeters have the form:
 \begin {gather*}
p(\cM,\pi_1)=10,\quad p(\cM,\pi_2)=12,\quad p(\cM,\pi_3)=12,\\
p(\cM,\pi_4)=13,\quad p(\cM,\pi_5)=13,\quad p(\cM,\pi_6)=12,
 \end {gather*}
hence, the maximal half-perimeter in this example is less than the weight of the minimal parametric filling again, and there is no an exact tour.
 \end {examp}

A minimal parametric filling $\cG=(G,\om)$ of a finite pseudo-metric space $\cM$ is said to be {\em unstable}, if there exists a minimal parametric filling $\cG'=(G',\om')$ of the same space $\cM$ such that $G'$ is obtained from $G$ by some splitting and $\om'(G')<\om(G)$. Otherwise, the filling $\cG$ is said to be
{\em stable}.

Further, notice that if a tree  $G'$ is obtained from the tree $G$ with a boundary $M$ by some splitting, then $\cO(G')\ss\cO(G)$, and for each proper splitting  $\cO(G')\ne\cO(G)$. Let $S$ be some set of splittings of the tree $G$ with the boundary $M$. For each $s\in S$, by $G_s$ we denote the result of the splitting $s$. We put $\cO_S(G)=\cup_{s\in S}\cO(G_s)$. If $\cG=(G,\om)$ is a minimal parametric filling and $S$ is the set of all its splittings decreasing the weight, then the set $\cO_S(G)$ is denoted by $\cO_d(\cG)$.

 \begin {ass}
If a minimal parametric filling $\cG$ of a finite pseudo-metric space $\cM$ is unstable, then any tour from $\cO_d(\cG)$ is not exact.
 \end {ass}

 \begin {proof}
Assume the contrary, i.e. let there exist an unstable minimal parametric filling $\cG$ of some pseudo-metric space $\cM=(M,\r)$ and some exact tour $\pi\in\cO_d(\cG)$. By $\cG'$ we denote a minimal parametric filling, whose weight is less than the weight of $\cG$, which is obtained from $\cG$ by some splitting, and such that
$\pi\in\cO(\cG')$.  Due to Assertion~\ref{ass:exact_order} and Corollary~\ref{cor:fill_ge_halfp}, $\om(\cG)=p(\cM,\pi)>\om(\cG')\ge
p(\cM,\pi)$, a contradiction. The proof is complete.
 \end {proof}

 \begin {examp}
The minimal parametric filling from the previous example is unstable. Its splitting into a minimal filling, the topology of the minimal filling together with the weight function are shown in Figure~\ref{fig:exam55}. Notice that this minimal filling possesses an exact tour, for example, the tour $(1,2,5,3,4)$.

\Pic{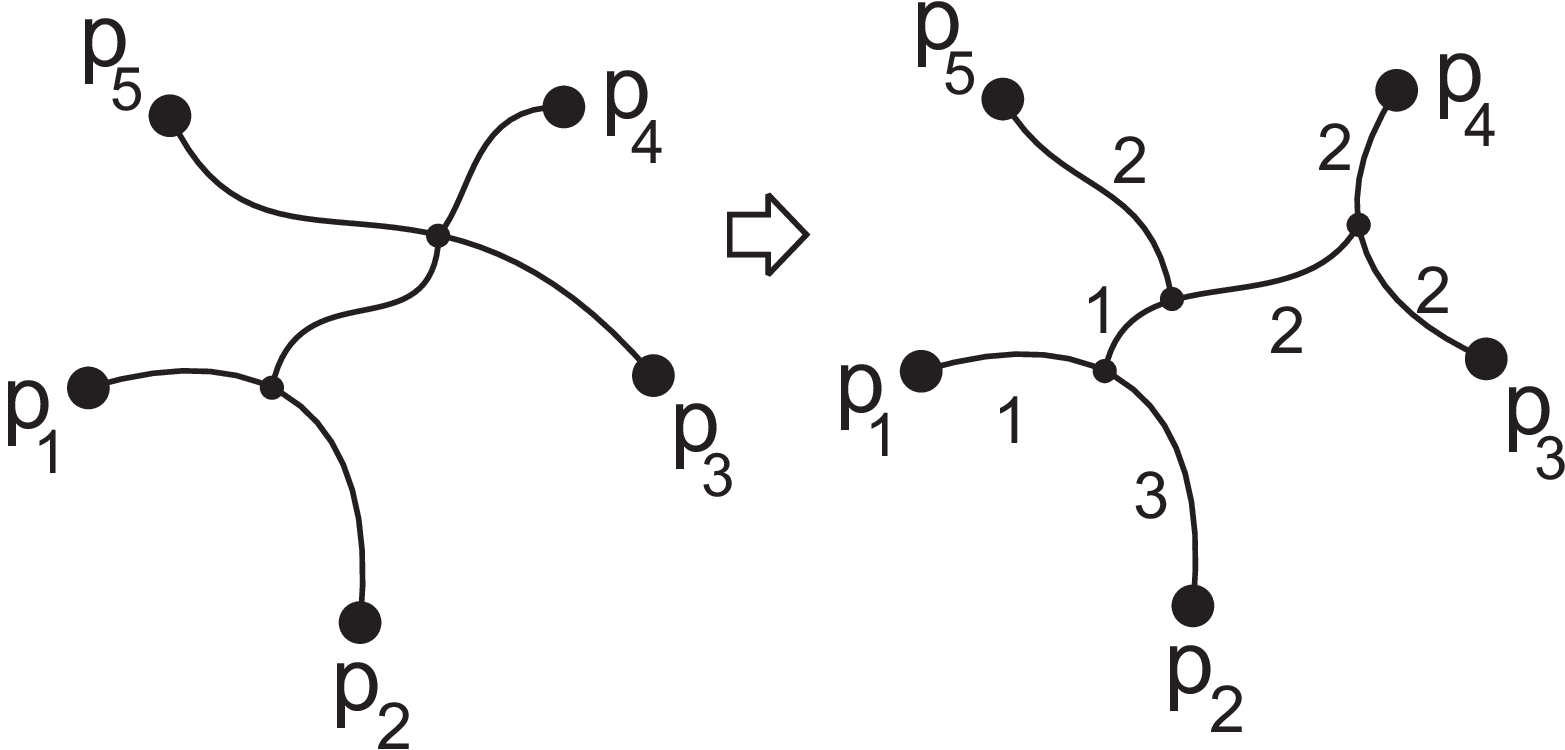}{Splitting of an unstable minimal parametric filling.}{fig:exam55}{200}
 \end {examp}

 \begin {conj}\label{conj:min-fill-exact-tour}
A minimal filling possesses an exact tour.
 \end {conj}

Let $\cG=(G,\om)$ be a weighted tree with boundary $M$. Cut $\cG$ over all boundary vertices of degree more than $1$. The resulting weighted trees  $\cG_i=(G_i,\om)$ with the boundaries $M_i=M\cap G_i$ are called by the {\em irreducible components\/} of the tree  $\cG$. We show that the presence of an exact tour is inherited by the irreducible components of the filling.

If $\cG=(G,\om)$ is a filling of a finite pseudo-metric space $\cM=(M,\r)$, then the irreducible components $\cG_i$ of the tree $\cG$ are also fillings of the pseudo-metric spaces $\cM_i=(M\cap G_i,\r)$, and if the filling $\cG$ is a minimal parametric one, then all $\cG_i$ are also minimal parametric fillings (notice, that the contrary is not true, generally speaking).

 \begin {ass}
Let $\cG=(G,\om)$ be a filling possessing an exact tour. Then each irreducible component $\cG_i$ of the tree $\cG$ also has an exact tour.
 \end {ass}

 \begin {proof}
We use induction over the number of the irreducible components. Let $\pi$ be an exact tour of the filling $\cG$ consisting of several irreducible components. By definition, the tour $\pi$ is generated by an embedding  $G'$ of the tree $G$. Let $y$ be some boundary vertex of the tree $G$, whose degree is more than $1$, and let $e_1,\ldots,e_m$, $m\ge2$, be the boundary edges incident to $y$. By $G_1,\ldots,G_m$ we denote the components obtained from $G$ by cutting over $y$. Each of them consists of less number of irreducible components than $G$ does. On the other hand, the vertices of the tree $G_i$ do not alternate with the vertices of the other components under the tour $\pi$ around the tree $G$. The latter means that, for some enumeration compatible with $\pi$, the boundary vertices of the component $G_i$, with the exception of $y$, have the form $(p_j,\ldots,p_{j+s})$. Therefore, all the paths joining consecutive vertices $p_k$, $p_{k+1}$, where $k=j,\ldots,(j+s-1)$, are also the paths from the tour $\pi$. Hence, they are exact. It remains to consider the paths joining $y$ and $p_j$, and $p_{j+s}$ and $y$ also. Each of those paths either is contained in the tour $\pi$ also, or is a sub-path of a boundary path from the tour $\pi$. In the both cases it is exact, thus, $G_i$ possesses an exact tour and, hence, due to the induction assumption, each irreducible component of the tree $G_i$ possesses an exact tour. It remains to notice that each irreducible component of the initial tree can be obtained in such a way. The proof is complete.
 \end {proof}

Let us state another conjecture and an important corollary.

 \begin {conj}\label{conj:min-fill-max-tour}
Let $\cM=(M,\r)$ be a pseudo-metric space joined by a tree $G$. Then
 $$
\mf(\cM)\le\max_{\pi\in\cO(G)}p(\cM,\pi).
 $$
 \end {conj}

Conjectures~\ref{conj:min-fill-exact-tour} and~\ref{conj:min-fill-max-tour} imply the following important result giving a formula to calculate the weight of a minimal filling of an arbitrary pseudo-metric space.

 \begin {conj}\label{conj:min-fill-formula}
For an arbitrary pseudo-metric space $\cM=(M,\r)$ the following formula is valid
 $$
\mf(\cM)=\min_G\max_{\pi\in\cO(G)}p(\cM,\pi),
 $$
where minimum can taken as over all trees with the boundary $M$, so as just over all binary trees with the boundary $M$.
 \end {conj}

 \section {Generalized Fillings and True Formula for the Weight of Minimal Filling} \label{sec:formula}
 \markright {\thesection.~Generalized Fillings and True Formula for the Weight of Minimal Filling.}
In this small section we give a review of the recent results obtained by our group when this text has been under preparation. We do not go in to details here referring to the original papers to appear soon, see~\cite{IOST} and~\cite{Eremin}.

 \subsection {Generalized fillings}
Investigating the fillings of pseudo-metric spaces, it turns out to be convenient to expand the class of weighted trees under consideration and to permit any weights of the edges (not only non-negative). The corresponding objects are called by {\em generalized fillings}, {\em minimal generalized fillings\/} and {\em minimal parametric ge\-ne\-ra\-lized fillings}. Their weights for a pseudo-metric space $\cM$ and a tree $G$ are denoted by $\mf_-(\cM)$ and $\mpf_-(\cM,G)$, respectively.

For any finite pseudo-metric space $\cM=(M,\r)$ and a tree $G$ joining $M$, the next evident inequality is valid: $\mpf_-(\cM,G)\le\mpf(\cM,G)$. And it is not difficult to construct an example, when this inequality becomes strict. It suffices to consider the four-points metric space consisting of the vertices of a unit square  and the binary tree $G$ whose mustaches join the diagonal vertices. If the negative weights are forbidden, then the interior edge has weight zero, see above. But in the case of generalized fillings one can decrease the total weight of the filling making the weight of the interior edge negative.

But for minimal generalized fillings the following result holds, see~\cite{IOST}.

 \begin {thm}\label{th:IOST}
Among minimal generalized fillings of an arbitrary finite pseudo-metric space $\cM$, there exists one with non-negative weight function, i.e. there exists a minimal filling. Therefore, $\mf_-(\cM)=\mf(\cM)$.
 \end {thm}

 \subsection {Multi-tours and the true formula for the weight of minimal filling}
Working with generalized fillings, A.~Eremin shows that the above Conjectures concerning exact tours are not true, but the formula from Conjecture~\ref{conj:min-fill-formula} can be modified to a valid one using the concept of multi-tour, introduced by him, see~\cite{Eremin}.

The concept of tour can be defined using Euler tour as follows. Let $G$ be a tree with a boundary $M$. Consider the doubling of $G$, i.e. the graph obtained from  $G$ by doubling of all its edges. This resulting graph possesses an Euler tour consisting of irreducible boundary paths. This Euler tour generates a cyclic order on $M$, which is a planar order generated by $G$. Evidently, each planar order can be obtained in such a way.

Now, instead of doubling $G$, let us consider the graph in which every edge of $G$ is taken with the multiplicity $2k$, $k\ge1$. The resulting graph possesses an Euler cycle consisting of irreducible boundary paths. This Euler cycle generates a {\em multi-tour of $M$ generated by $G$}. Such multi-tour can be considered as a bijection $\pi\:X\to X$, where $X=\sqcup_{i=1}^kM$. The set of all multi-tours generated by $G$ on $M$ is denoted by $\cO_\mu(M,G)$.

Let $\cM=(M,\r)$ be a finite pseudo-metric space, and $G$ be a tree joining $M$. As in the case of tours, let us define the {\em half-perimeter of $M$ corresponding to a multi-tour $\pi$ generated by $G$} putting
 $$
p(\cM,G,\pi)=\frac1{2k}\sum_{x\in X}\r\bigl(x,\pi(x)\bigr).
 $$

Eremin obtained the following formula calculating the weight of minimal parametric generalized filling and minimal filling.

 \begin {thm}[A.~Eremin]\label{th:eremin}
For an arbitrary finite pseudo-metric space $\cM=(M,\r)$ and an arbitrary tree $G$ joining $M$, the weight of minimal parametric generalized filling can be calculated as follows
 $$
\mpf_-(\cM,G)=\max\bigl\{p(\cM,G,\pi)\mid \pi\in\cO_\mu(M,G)\bigr\}.
 $$
The weight of minimal filling can be calculated as follows
 $$
\mf(\cM)=\mf_-(\cM)=\min_G\max\bigl\{p(\cM,G,\pi)\mid \pi\in\cO_\mu(M,G)\bigr\},
 $$
where minimum is taken over all binary trees $G$ joining $M$.
 \end{thm}

Theorem~\ref{th:eremin} gives an opportunity to prove several interesting corollaries. For example, it is shown in~\cite{Eremin} that minimal filling of a finite pseudo-metric space {\sl in general position\/} is a binary tree with positive weights.

 \section {Additive Spaces} \label{sec:additive}
 \markright {\thesection.~Additive Spaces.}
The additive spaces are very popular in bioinformatics, playing an important role in evolution theory. Recall that a finite metric space $\cM=(M,\r)$ is called {\em additive}, if $M$ can be joined by a weighted tree $\cG=(G,\om)$ such that $\r$ coincides with the restriction of $d_\om$ onto $M$. The tree $\cG$ in this case is called {\em generating tree\/} for the space $\cM$.

Not any pseudo-metric space is additive. It turns out that an additivity criterion can be stated in terms of well-known {\em $4$ points rule}: for any four points $p_i$, $p_j$, $p_k$, $p_l$, the values $\r(p_i,p_j)+\r(p_k,p_l)$, $\r(p_i,p_k)+\r(p_j,p_l)$, $\r(p_i,p_l)+\r(p_j,p_k)$ are the lengths of sides of an isosceles triangle whose base does not exceed its other sides.

 \begin {ass}{\rom{(\cite{Zaretskij}, \cite{SimoesPereira})}}
A pseudo-metric space is additive, if and only if it meets the $4$ points rule.
 \end {ass}

The next statement deals with the uniqueness question for the generating tree.

 \begin {ass}{\rom{(\cite{Smolenskij}, \cite{HakimiYau})}}\label{ass:add-unique}
In the class of non-degenerate weighted trees, the generating tree of an additive metric space is unique.
 \end {ass}

The next criterion solves completely the minimal filling problem for additive pseudo-metric spaces.

 \begin {thm}\label{th:additive=minimum}
Minimal fillings of an additive pseudo-metric space are exactly its generating trees.
 \end {thm}

 \begin {proof}
At first, let $\cG=(G,\om)$ be a generating tree of a finite additive pseudo-metric space $\cM=(M,\r)$. We show that $\cG$ is a minimal filling for $\cM$.

Let $\cG'=(G',\om')$ be some minimal filling of $\cM$, and $\pi'$ be a planar order on $M$ with respect to $G'$. Then, due to Assertion~\ref{ass:2weights},
 $$
\sum_{p\in M}d_{\om'}\bigl(p,\pi'(p)\bigr)=2\om'(G').
 $$
On the other hand, for any $p$ and $q$ from $M$, we have $d_{\om'}(p,q)\ge\r(p,q)=d_\om(p,q)$, where the inequality follows from the definition of a filling, and the equality is valid since $\cG$ is a generating tree of the additive space $\cM$. Therefore,
 $$
\sum_{p\in M}d_{\om}\bigl(p,\pi'(p)\bigr)\le
\sum_{p\in M}d_{\om'}\bigl(p,\pi'(p)\bigr)=
2\om'(G').
 $$
Again due to Assertion~\ref{ass:2weights}, we have
 $$
2\om(G)\le\sum_{p\in M}d_{\om}\bigl(p,\pi'(p)\bigr),
 $$
and hence $\om(G)\le\om(G')$, and so $\cG$ is also a minimal filling of the space $\cM$, q.e.d.

Now let us prove the inverse statement. Let $\cM=(M,\r)$ be an additive space, and $\cG=(G,\om)$ be its arbitrary minimal filling of $\cM$. We show that $\cG$ is a generating tree for  $\cM$.

Let $\ti\cG=(\ti G,\ti\om)$ be some generating tree for $\cM$. Due to the direct statement of Theorem already proved, the tree $\ti\cG$ is a minimal filling for $\cM$, and hence $\om(G)=\ti\om(\ti G)$. Let $\pi$ be some tour around the tree $G$. Then, due to Assertion~\ref{ass:2weights},
 $$
2\ti\om(\ti G)\le\sum_{p\in M}d_{\ti\om}\bigl(p,\pi(p)\bigr)
\qquad\text{и}\qquad
\sum_{p\in M}d_\om\bigl(p,\pi(p)\bigr)=2\om(G),
 $$
and therefore, due to the condition $\r\bigl(p,\pi(p)\bigr)=d_{\ti\om}\bigl(p,\pi(p)\bigr)\le d_\om\bigl(p,\pi(p)\bigr)$ and the equality  $\om(G)=\ti\om(\ti G)$, we conclude, that $\r\bigl(p,\pi(p)\bigr)=d_\om\bigl(p,\pi(p)\bigr)$ for any point $p\in M$. On the other hand, for any pair $p$ and $q$ of the vertices $M$, there exists a tour of the tree $G$, such that these points are consecutive in it, and hence, $\r(p,q)=d_\om(p,q)$ for any $p$ and $q$ from $M$, i.e\. the tree $\cG$ is also generating for $\cM$. Theorem is proved.
 \end {proof}

Let $\cG=(G,\om)$ be an arbitrary weighted tree with a boundary $M$. It is evident, then the tree $\cG$ is a filling of the pseudo-metric space $(M,d_\om)$, and, hence, due to Theorem~\ref{th:additive=minimum}, it is a minimal filling of this space. On the other hand, each minimal filling has such a form in accordance with the following result.

 \begin {ass}\label{ass:induced_minimum}
Let $\cG=(G,\om)$ be a minimal filling of a finite pseudo-metric space $\cM=(M,\r)$. Then $\cG$ is a minimal filling of the space $(M,d_\om)$ also.
 \end {ass}

 \begin {proof}
Indeed, $\cG$ is a generating tree of the additive space $(M,d_\om)$. It remains to apply Theorem~\ref{th:additive=minimum}. The proof is complete.
 \end {proof}

 \begin {cor}
The set of all minimal fillings coincides with the set of all weighted trees.
 \end {cor}

 \begin {cor}
Let $\G$ be some network parameterized by a tree $G$ and joining a finite subset $M$ of a pseudo-metric space $\cX$. By $\om_\G$ we denote the weight function on $G$, generated by the network $\G$. Then $(G,\om_\G)$ is a minimal filling of the space $(M,d_{\om_\G})$.
 \end {cor}

Theorem~\ref{th:additive=minimum} can be reformulated as the following description of generating trees of an additive pseudo-metric space.

 \begin {cor}\label{cor:fill=gen}
A weighted tree $\cG$ joining an additive space $\cM$ is a generating tree for $\cM$, if and only if $\cG$ is a minimal filling of the space $\cM$.
 \end {cor}

Recall that boundary rigidity problem for a metric in Riemannian geometry is a uniqueness problem for the reconstruction of this metric from the interior metric induced on the boundary of the initial manifold. In the scope of this paper, the boundary rigidity problem is a uniqueness problem for non-degenerate minimal filling for a given finite pseudo-metric space. More rigorous, a pseudo-metric $\r$ on a finite set $M$ is said to be {\em boundary rigid}, if minimal filling of $(M,\r)$ is unique in the class of non-degenerate trees.

Assertion~\ref{ass:add-unique} and Corollary~\ref{cor:fill=gen} imply the following result directly.

 \begin {cor}\label{cor:rigid}
Each additive metric $\r$ on a finite space $M$ is boundary rigid.
 \end {cor}

We give another couple of results concerning additive spaces and their minimal fillings.

 \begin {cor}\label{cor:weight-mf-additive}
The weight of a minimal filling of an additive space $\cM=(M,\r)$ is equal to the half-perimeter of the space $\cM$.
 \end {cor}

 \begin {proof}
Let $\pi$ be an arbitrary tour of a generating tree $(G,\om)$ of the space $\cM$, and $\s$ be an arbitrary cyclic order on $M$. Since $\cM$ is additive, then the tour $\pi$ is exact, therefore $\om(G)=p(\cM,\pi)$ in accordance with Assertion~\ref{ass:exact_order}. Since $\r=d_\om$, Assertion~\ref{ass:2weights} implies that $\om(G)\le
p(\cM,\s)$. Therefore, $p(\cM,\pi)=\min_\s p(\cM,\s)$, where minimum is taken over all cyclic orders $\s$ on $M$. Thus, $p(\cM,\pi)=p(\cM)$. It remains to apply Theorem~\ref{th:additive=minimum} that implies $\om(G)=\mf(\cM)$. The proof is complete.
 \end {proof}

The next additivity criterion is obtained by O.~Rubleva, a student of mechanical and mathematical faculty of Moscow state university, see~\cite{Rubleva}.

 \begin {ass}[O.~Rubleva]\label{ass:Rubleva}
The weight of a minimal filling of a finite pseudo-metric space is equal to the half-perimeter of this space, if and only if this space is additive.
 \end {ass}

Let $(M,\r)$ be an arbitrary finite pseudo-metric space. We say that its minimal filling $(G,\om)$ is {\em inherited}, if for any non-empty subset $S\ss M$ the subtree $(G_S,\om)$ in $(G,\om)$ consisting of all the paths in $G$ joining the vertices from $S$ is a minimal filling of the pseudo-metric space $(S,\r)$, where the restriction of the metric $\r$ onto $S$ and the restriction of the weight function $\om$ onto $G_S$ are denoted by the same letters.

 \begin {cor}\label{cor:add_restrict}
If a finite pseudo-metric space $\cM$ is additive, then any its minimal filling is inherited. Conversely, if some minimal filling of the space $\cM$ is inherited, then the space $\cM$ is additive and, hence, any its minimal filling is inherited.
 \end {cor}

 \begin {proof}
Let $\cM=(M,\r)$ be additive and $S$ be an arbitrary non-empty subset of $M$. In this case, due to Theorem~\ref{th:additive=minimum}, any minimal filling $(G,\om)$ of the space $\cM$ is a generating tree, therefore, for any points $x$ and $y$ from $S$ and the path $\g$ joining them in $G_S$, the equality $\r(x,y)=\om(\g)$ is valid, and so $(S,\r)$ is additive space with generating tree $(G_S,\om)$.  In accordance with Theorem~\ref{th:additive=minimum}, this tree is a minimal filling of $(S,\r)$.

Conversely, let some minimal filling $\cG=(G,\om)$ of the space $\cM$ is inherited. Then, for any two points $x$ and $y$ from $M$, consider the set $S=\{x,y\}$. The subtree $G_S$ for $S$ coincides with the path $\g_{xy}$ joining $x$ and $y$ in $G$ and, due to our assumptions, is a minimal filling of the space  $(S,\r)$, so the weight of the path $\g_{xy}$ is equal to $\r(x,y)$. The latter implies, that the tree $\cG$ is generating for the space $\cM$, i.e. $\cM$ is additive. Corollary is proved.
 \end {proof}

\begin {rk}
In the scope of Assertion~\ref{ass:Rubleva}, we conjectured that if there exists a tree joining a pseudo-metric space such that all the corresponding pereimeters are equal to each other, then the space is additive. It turns out that it is not true. Z.~Ovsjannikov suggested to consider a wider class of spaces, so called pseudo-additive spaces, for which our conjecture becomes true. Here we just announce this results briefly, referring to paper~\cite{Ovs}.

A finite pseudo-metric space $\cM=(M,\r)$ is said to be {\em pseudo-additive}, if the metric $\r$ coincides with $d_\om$ for a weighted tree $(G,\om)$ (which is aslo called {\em generating\/}), where the weight function $\om$ can take arbitrary (not necessary non-negative) values. Z.~Ovsjannikov shows that these spaces can be described in terms of so-called {\em weak $4$-points rule}: for any four points $p_i$, $p_j$, $p_k$, $p_l$, the values $\r(p_i,p_j)+\r(p_k,p_l)$, $\r(p_i,p_k)+\r(p_j,p_l)$, $\r(p_i,p_l)+\r(p_j,p_k)$ are the lengths of sides of an isosceles triangle. The generating tree is also unique in the class of non-degenerate trees.  Moreover, the following result is valid, see~\cite{Ovs}.

 \begin{thm}[Z.~Ovsjannikov]\label{th:ovs}
Let $\cM=(M,\r)$ be a finite pseudo-metric space. Then the following statements are equivalent.
 \begin{itemize}
 \item There exist a tree $G$ such that $M$ coincides with the set of degree $1$ vertices of $G$ and all the half-perimeters $p(M,G,\pi)$ of $M$ corresponding to the tours of $G$ are equal to each other.
 \item The space $\cM$ is pseudo-additive.
 \end{itemize}
Moreover, the three $G$ in this case is a generating tree for the space $\cM$.
 \end{thm}

It would be interesting to see what role could play these pseudo-additive spaces in applications.
 \end {rk}

 \section {Rays of Pseudo-Metric Spaces and Minimal Fillings} \label{sec:rays}
 \markright {\thesection.~Rays of Pseudo-Metric Spaces and Minimal Fillings.}
Let $\cM=(M,\r)$ be a pseudo-metric space and $\l$ be a positive real number. By $\l\r$ we denote the function on the pairs of points from $M$, defined as follows:  $(\l\r)(x,y)=\l\r(x,y)$. It is clear, that $\l\r$ is a pseudo-metric $\l$, and all the spaces $(M,\l\r)$ are degenerate or non-degenerate simultaneously. It is also clear, that $(G,\om)$ is a minimal (parametric) filling of the space $(M,\r)$, if and only if $(G,\l\om)$ is a minimal (parametric) filling of the space $(M,\l\r)$. The set of all pseudo-metric spaces $(M,\l\r)$, $\l>0$, we call by the {\em open multiplicative ray\/} passing through $\cM$. The {\em closed multiplicative ray\/} is the set of the spaces $(M,\l\r)$, $\l\ge0$. It can be obtained from the open one by adding the trivial pseudo-metric space which all the distances are zero in.

By $\r+a$ we denote the function on the pairs of points from $M$ defined as follows: $(\r+a)(x,y)=\r(x,y)+a$.  It is easy to see, that there exists a value $a_\cM\le0$ such that the function $\r+a$ is a pseudo-metric for $a\ge a_\cM$, and is not a pseudo-metric for $a<a_\cM$. The set of all the pseudo-metric spaces $\cM+a$, $a\ge a_\cM$, we call the {\em additive ray with the vertex at $\cM+a_\cM$, passing through $\cM$} and denote by $[\cM)$. An additive ray without its vertex is called  {\em open}. Notice that a pseudo-metric space is not a vertex of an additive ray, if and only if it is non-degenerate. In particular, all the spaces from an arbitrary open additive ray are non-degenerate.

Let $\cG=(G,\om)$ be a weighted tree joining the points of a pseudo-metric space $\cM=(M,\r)$, and $M$ coincides with the set of degree $1$ vertices of the tree $G$ and consists of $n$ points. By $b$ we denote the smallest weight of the boundary edges of the tree $\cG$. Let $a\ge a_\cM$ be such that the inequality $b+a/2\ge0$ holds. Define a weight function $\om_a$ on $G$ adding the value $a/2$ to the weights of all the boundary edges of the tree $\cG$, and let $\cG_a=(G,\om_a)$.

 \begin {ass}\label{ass:ray-a-negative}
Under the above notations, assume in addition that $a\le0$ and $\cG$ is a minimal filling of the space $\cM$. Then $\cG_a$ is a minimal filling of the space $\cM+a$.
 \end {ass}

 \begin {proof}
It is easy to see that $\cG$ and $\cG_a$ are ore are not some fillings of the corresponding spaces  $\cM$ and $\cM+a$ simultaneously. Besides, $\om_a(G)=\om(G)+an/2$.

Assume the contrary, i.e. let the filling  $\cG_a$ be not minimal. Consider a minimal filling  $\cG'=(G',\om')$ of the space $\cM+a$, which is a binary tree (it exists due to Theorem~\ref{th:trees_fill}). Then $\om'(G')<\om_a(G)=\om(G)+an/2$. Put $c=-a$. Since $a\le0$, then $c\ge0$, and hence, if we add $c/2$ to the $\om'$-weights of the boundary edges of the tree $G'$, then we obtain a weight function $\om'_c$ on $G'$, and $\cG'_c=(G',\om'_c)$ is a filling of the space $\cM$. But the weight of the filling obtained is equal to $\om'(G')+cn/2<\om(G)+an/2+cn/2=\om(G)$, that contradicts to the minimality of the filling $\cG$. The proof is complete.
 \end {proof}

 \begin {ass}\label{ass:ray-inequality}
Under the above assumptions, assume in addition that  $\cG$ is a minimal filling of the space $\cM$. Then $-2b\le a_\cM$.
 \end {ass}

 \begin {proof}
Assume the contrary, i.e\.  $-2b>a_\cM$. Then the function $\r-2b$ is a non-degenerate metric on $M$, i.e\. $\cM-2b$ is a non-degenerate metric space. Since $-2b\le0$, then, due to Assertion~\ref{ass:ray-a-negative}, the weighted tree $\cG_{-2b}=(G,\om_{-2b})$ is a minimal filling of the non-degenerate space $\cM-2b$.  But the tree $\cG_{-2b}$ contains a degenerate boundary edge, that contradicts to Corollary~\ref{cor:boundary-vertices-nondegenerate}. The proof is completed.
 \end {proof}

 \begin {ass}
Under the above notations, assume in addition that  $a>0$ and $\cG$ is a minimal filling of the space  $\cM$. Then $\cG_a$ is a minimal filling of the space $\cM+a$.
 \end {ass}

 \begin {proof}
Assume the contrary, i.e\. let the filling $\cG_a$ be not minimal, and let $\cG'=(G',\om')$ be a minimal filling of the space $\cM+a$, which is a binary tree (it exists due to Theorem~\ref{th:trees_fill}). Then $\om'(G')<\om_a(G)=\om(G)+an/2$. By $b'$ we denote the minimal weight of the boundary edges of the tree $\cG'$. Due to Assertion~\ref{ass:ray-inequality}, we have $-2b'\le a_\cM+a$, therefore $-2b'\le a$ and the function $\om'_{-a}$ is a weight function on $G'$, and the weighted tree $\cG'_{-a}$ is a filling of the initial space $\cM$.  But then $\om'_{-a}(G')=\om'(G')-an/2<\om(G)+an/2-an/2=\om(G)$, that contradicts to minimality of the filling $\cG$. The proof is complete.
 \end {proof}

Let us summarize the results of this section in the following Theorem.

 \begin {thm}\label{th:rays}
Let $\cG=(G,\om)$ be a minimal filling of a pseudo-metric space $\cM$, and $M$ coincides with the set of degree $1$ vertices of the tree $G$. Let $b$ be the least weight of the boundary edges of the filling $\cG$. Then $-2b<a_\cM$, and hence, for any $a\ge a_\cM$, the function $\om_a$ obtained from $\om$ by adding $a/2$ to all the weights of the boundary edges is a weight function. Moreover, for any such $a$, the weighted tree $\cG_a=(G,\om_a)$ is a minimal filling of the space $\cM+a$. In particular, all the spaces of the form $(M,\l\r+a)$, $\l>0$, $a>\l a_\cM$, have the same set of types of minimal fillings, and all these types are realized among the types of minimal fillings of the spaces $(M,\l\r+\l a_\cM)$, $\l\ge0$.
 \end {thm}

 \section {Incommensurable Spaces} \label{sec:incommens}
 \markright {\thesection.~Incommensurable Spaces.}
A finite pseudo-metric space in which all non-zero distances are linear independent over the field $\Q$ of rational numbers is said to be {\em incommensurable}. Notice that each incommensurable metric space is non-degenerate.

 \begin {ass}\label{ass:noncomp}
An additive space consisting of $n>3$ points can not be incommensurable.
 \end {ass}

 \begin {proof}
The distances between points of an additive space lies in the vector space over $\Q$, spanned onto the weights of the edges of the generating tree. The dimension of this space is at most  $2n-3$. On the other hand, if this space is incommensurable, then the dimension of the linear hull over $\Q$ of all the distances is equal to  $n(n-1)/2$ which is more than $2n-3$ for $n>3$. The proof is complete.
 \end {proof}

 \begin {cor}\label{cor:add-noncomparable}
A minimal filling of an incommensurable metric space consisting of more than three points contains a non-degenerate interior edge.
 \end {cor}

 \begin {proof}
If it is not so, then, due to non-degeneracy of the space and Corollary~\ref{cor:boundary-vertices-nondegenerate}, the minimal filling is a star of the unique interior vertex and, hence, due to Corollary~\ref{cor:star}, the space is additive, that contradicts to Assertion~\ref{ass:noncomp}. Corollary is proved.
 \end {proof}

 \section {Minimal Fillings Examples} \label{sec:examp}
 \markright {\thesection.~Minimal Fillings Examples.}
In this section we give several examples of minimal filling and demonstrate, how to use the technique elaborated above.

 \subsection {Triangle}\label{subsec:triangle}
Let $\cM=(M,\r)$ consist of three points $p_1$, $p_2$, and $p_3$. Put $\r_{ij}=\r(p_i,p_j)$. Consider the tree $G=(V,E)$ with $V=M\cup\{v\}$ and $E=\{vp_i\}_{i=1}^3$. Define the weight function $\om$ on $E$ by the following formula:
 $$
\om(e_i)=\dfrac{\r_{ij}+\r_{ik}-\r_{jk}}2,
 $$
where $\{i,j,k\}=\{1,2,3\}$. Notice that $d_\om$ restricted onto $M$ coincides with $\r$. Therefore, $\cM$ is an additive space, $\cG=(G,\om)$ is a generating tree for $\cM$, and, due to Theorem~\ref{th:additive=minimum}, $\cG$ is a minimal filling of $\cM$.

Notice that in this case $\om(G)=p(\cM)$, each tour is exact, and the tree $\cG$ is non-degenerate if and only if $\cM$ is a non-degenerate triangle.

Recall that the value $(\r_{ij}+\r_{ik}-\r_{jk})/2$ is called by the {\em Gromov product\/} $(p_j,p_k)_{p_i}$ of the points $p_j$ and $p_k$ of the space $\cM$ with respect to the point $p_i$, see~\cite{GromHypGr}. The next result is useful in particular calculations.

 \begin {ass}
The identity
 $$
\r_{ik}=(p_j,p_k)_{p_i}+(p_i,p_j)_{p_k}
 $$
is valid.
 \end {ass}

In what follows we also need the following technical result.

 \begin {ass}\label{ass:seim-mf-3}
Let $\cG=(G,\om)$ be a filling of a pseudo-metric space $\cM=(M,\r)$ consisting of at least $3$ points. Let $e=p_1v$ be an edge from $G$, and $p_1\in M$, $v\not\in M$. Consider two some boundary paths $\g_{1i}$, $i=2,\,3$, starting at $p_1$ and containing $e$. By $p_i$ we denote the boundary vertex, which the path
$\g_{1i}$ comes to. By $\g_{23}$ we denote the path joining $p_2$ and $p_3$, and by $\r_{ij}$ --- the distance between $p_i$ and $p_j$. Then
 \begin {itemize}
  \item if the paths $\g_{12}$ and $\g_{13}$ are exact, then $\om(e)\le(p_2,p_3)_{p_1}$\/\rom;
  \item if the path $\g_{23}$ and one of the paths $\g_{1i}$ are exact and the paths $\g_{12}$ and $\g_{13}$ have exactly one common edge, then $\om(e)\ge(p_2,p_3)_{p_1}$.
 \end {itemize}
 \end {ass}

 \begin {proof}
In the first case, since the edge $e$ is contained in the intersection of the paths $\g_{12}$ and $\g_{13}$, and the path $\g_{23}$ is equal to their symmetric difference, then
 $$
\om(e)\le\frac{\om(\g_{12})+\om(\g_{13})-\om(\g_{23})}{2}=
\frac{\r_{12}+\r_{13}-\om(\g_{23})}{2}\le(p_2,p_3)_{p_1},
 $$
where the equality is valid due to our assumptions concerning exactness, and the second inequality follows from the definition of filling.

In the second case, the edge $e$ is equal to the intersection of the paths  $\g_{12}$ and $\g_{13}$, so
 $$
\om(e)=\frac{\om(\g_{12})+\om(\g_{13})-\om(\g_{23})}{2}=
\frac{\r_{1i}+\om(\g_{1j})-\r_{23}}{2}\ge(p_2,p_3)_{p_1},
 $$
where $\{i,j\}=\{2,3\}$, the second equality is valid due to our assumptions concerning exactness, and the inequality is true to the definition of filling. Assertion is proved.
 \end {proof}

 \subsection {Regular Simplex}\label{subsec:simplex}
Let all the distances in the metric space $\cM$ are the same and are equal to $d$, i.e\.  $\cM$ is a regular simplex. Then the weighted tree $\cG=(G,\om)$, $G=(V,E)$, with the vertex set $V=M\cup\{v\}$ and edges $vm$, $m\in M$, of the weight $d/2$ is generating for $\cM$. Therefore, the space $\cM$ is additive, and, due to Theorem~\ref{th:additive=minimum}, $\cG$ is its unique non-degenerate minimal filling. If $n$ is the number of points in $M$, then the weight of the minimal filling is equal to $dn/2$.

 \subsection {Star}
If a minimal filling $\cG=(G,\om)$ of a space $\cM=(M,\r)$ is a star whose single interior vertex $v$ is joined with all the points $p_i\in M$, $1\le i\le n$, $n\ge 3$, then, due to Corollary~\ref{cor:star}, the pseudo-metric space $\cM$ is additive. In this case the weights of edges can be calculated easily. Indeed, put $e_i=vp_i$. Since a subspace of an additive space is additive itself, then we can use the results for three-points additive space, see above. So, we have $\om(e_i)=(p_j,p_k)_{p_i}$, where $i$, $j$, and $k$ are arbitrary distinct numbers of the boundary vertices, and $\r_{ij}$ is the distance between the corresponding points.

Notice that in accordance with Corollary~\ref{cor:weight-mf-additive} in this case the equality $\om(G)=p(\cM)$ is valid, any tour is exact, and the tree  $\cG$ is non-degenerate if and only if $\cM$ is a non-degenerate space.

 \subsection {Mustaches of Degree more than Two}
In accordance with Corollary~\ref{cor:moustaches-additive-subspace}, any mustaches of a minimal filling of a pseudo-metric space forms an additive subspace. If the degree of such mustaches is more than two, then we can calculate the weights of all the edges containing in the mustaches just in the same way as in the case of a star. Notice that the edges of the considered mustaches are non-degenerate, if and only if the additive subspace genet=rated by the mustaches is non-degenerate.

 \subsection {Two Stars}
Now assume that a minimal filling $\cG$ of a space $\cM$ possesses two adjacent interior vertices  $v_1$ and $v_2$, each of which is adjacent with three or more vertices of degree $1$. Then the weights of the boundary vertices can be calculated as described above. By $\om_{ij}$ we denote the weight of the edge joining the $i$th boundary vertex $p_i$ with  $v_j$. Then the weight of the edge $v_1v_2$ is the least among the values meeting all the inequalities of the form $\om_{i1}+\om(v_1v_2)+\om_{j2}\ge\r_{ij}$, where $\r_{ij}=\r(p_i,p_j)$ as above. Therefore, $\om(v_1v_2)=\max_{i,j}(\r_{ij}-\om_{i1}-\om_{j2})$, where maximum is taken over all the pairs $(p_i,p_j)$ of boundary vertices such that $p_i$ is adjacent with $v_1$ and $p_j$ is adjacent with $v_2$.

 \subsection {Four-Points Spaces}
Here we give a complete description of minimal fillings for four-points spaces. We start with the following result, which follows from Corollaries~\ref{cor:any_pair} and~\ref{cor:4-configuration-mf} immediately.

 \begin {cor}
For four-points pseudo-metric spaces, each minimal filling possesses an exact tour.
 \end {cor}

 \begin {ass}
Let $M=\{p_1,p_2,p_3,p_4\}$, and $\r$ be an arbitrary pseudo-metric on $M$. Put $\r_{ij}=\r(p_i,p_j)$. Then the weight of a minimal filling $\cG=(G,\om)$ of the space $\cM=(M,\r)$ is given by the following formula
 $$
\frac12\bigl(\min\{\r_{12}+\r_{34},\,\r_{13}+\r_{24},\,\r_{14}+\r_{23}\}+
\max\{\r_{12}+\r_{34},\,\r_{13}+\r_{24},\,\r_{14}+\r_{23}\}\bigr).
 $$
If the minimum in this formula is equal to $\r_{ij}+\r_{rs}$, then the type of minimal filling is the binary tree with the mustaches $\{p_i,p_j\}$ and $\{p_r,p_s\}$.
 \end {ass}

 \begin {proof}
Let $\cG=(G,\om)$ be a minimal filling of the space $\cM$, which is a binary tree. By $\l_i$ we denote the weight of the edge incident to $p_i$, and by $\l_5$ --- the weight of the remaining edge of the tree $\cG$. Without loss of generality, assume that $\{p_1,p_2\}$ are mustaches of the tree $G$.

In accordance with Corollary~\ref{cor:any_pair}, we have $\l_1+\l_2=\r_{12}$ and $\l_3+\l_4=\r_{34}$.

Due to Corollary~\ref{cor:4-configuration-mf}, for some $i$ and $j$, $\{i,j\}=\{3,4\}$, the points $p_1$ and $p_i$, and the points $p_2$ and $p_j$ also are joined by exact paths. Without loss of generality assume that $p_1$ and $p_3$, and $p_2$ and $p_4$ are joined by the exact paths. Then $\l_1+\l_3+\l_5=\r_{13}$ and $\l_2+\l_4+\l_5=\r_{24}$, and so, taking the alternate sum of the four equalities, where the latter two are taken with plus, and the first two are taken with minus, we conclude that
 $$
\l_5=\dfrac{\r_{13}+\r_{24}-\r_{12}-\r_{34}}2.
 $$
In particular, since $\l_5\ge0$, we have $\r_{13}+\r_{24}\ge\r_{12}+\r_{34}$. Besides,
 $$
\om(G)=\dfrac{\r_{13}+\r_{24}+\r_{12}+\r_{34}}2.
 $$

Further, $\l_1+\l_5+\l_4\ge\r_{14}$ and $\l_2+\l_5+\l_3\ge\r_{23}$, and so, summing these inequalities and substituting the expression for $\l_5$, we obtain
 $$
\l_1+\l_4+\l_5+\l_2+\l_3+\l_5=
\r_{12}+\r_{34}+(\r_{13}+\r_{24}-\r_{12}-\r_{34})\ge\r_{14}+\r_{23},
 $$
i.e\. $\r_{13}+\r_{24}=\max\{\r_{12}+\r_{34},\,\r_{13}+\r_{24},\,\r_{14}+\r_{23}\}$. The latter condition we call by Condition~$(*)$.

We put $\l_1=x$. Then, due to Assertion~\ref{ass:seim-mf-3}, we have
 $$
(p_2,p_4)_{p_1}\le x\le(p_2,p_3)_{p_1}.
 $$
Thus, under the above assumptions, the weight function of the minimal filling  $\cG$ meets the following equations system, which we call by System~$(**)$:
 \begin {gather*}
(p_2,p_4)_{p_1}\le\l_1=x\le(p_2,p_3)_{p_1},\ \
\l_2=\r_{12}-x,\\
\l_3=\dfrac{\r_{13}-\r_{24}+\r_{12}+\r_{34}}2-x,\\
\l_4=\dfrac{-\r_{13}+\r_{24}-\r_{12}+\r_{34}}2+x,\\
\l_5=\dfrac{\r_{13}+\r_{24}-\r_{12}-\r_{34}}2.
 \end {gather*}
Show that each solution of this system gives some filling of the space $\cM$, providing Condition~$(*)$ is valid.To do this it suffices to verify that $(p_2,p_4)_{p_1}\le(p_2,p_3)_{p_1}$, all $\l_i$, $1\le i\le4$, are non-negative, and also that $\l_1+\l_5+\l_4\ge\r_{14}$ and $\l_2+\l_5+\l_3\ge\r_{23}$.

The inequality  $(p_2,p_4)_{p_1}\le(p_2,p_3)_{p_1}$ follows from Condition~$(*)$. Further,
 \begin {gather*}
\l_1\ge (p_2,p_4)_{p_1}\ge0;\\
\l_2\ge\r_{12}-(p_2,p_3)_{p_1}=(p_1,p_3)_{p_2}\ge0;\\
\l_3\ge\dfrac{\r_{13}-\r_{24}+\r_{12}+\r_{34}}2-(p_2,p_3)_{p_1}=
(p_3,p_4)_{p_2}\ge0;\\
\l_4\ge\dfrac{-\r_{13}+\r_{24}-\r_{12}+\r_{34}}2+(p_2,p_4)_{p_1}=
(p_1,p_3)_{p_4}\ge0;\\
\l_1+\l_5+\l_4=\r_{24}-\r_{12}+2x\ge\r_{24}-\r_{12}+2(p_2,p_4)_{p_1}
=\r_{14};\\
\l_2+\l_5+\l_3=\r_{12}+\r_{13}-2x\ge\r_{12}+\r_{13}-2(p_2,p_3)_{p_1}
=\r_{23},
 \end {gather*}
q.e.d.

To complete the proof of Assertion it remains to explore the information obtained.

 \begin {lem}
Let $\cM=(M,\r)$ be a four-point pseudo-metric space such that Condition~$(*)$ is valid, and let $G$ be a binary tree with the boundary $M$, and ${p_1,\,p_2}$ are mustaches of the tree $G$. Define the numbers $\l_i$ from System~$(**)$.  Then all $\l_i$ are non-negative and, hence, determine a weight function $\om$, and the weighted tree $\cG=(G,\om)$ is a minimal parametric filling of the weight $(\r_{12}+\r_{24}+\r_{43}+\r_{31})/2$.
 \end {lem}

 \begin {proof}
The verification of non-negativity of all $\l_i$ can be performed just in the same way as in the investigation of a minimal filling. Further, similarly one can show that $\cG$ is a filling of the space $\cM$, and the tour with the enumeration $(p_1,\,p_2,\,p_4,\,p_3)$ is exact. So, due to Assertion~\ref{ass:exact_order}, the tree $\cG$ is a minimal parametric filling. The weight of this filling can be calculated directly. Lemma is proved.
 \end {proof}

Similar formulas take place under the assumptions that the pair $\{p_2,\,p_3\}$ is a mustaches of $G$ and that $\cM$ meets Condition~$(*)$. In this case we aslo obtain a minimal parametric filling. The weight of this filling is equal to $(\r_{23}+\r_{31}+\r_{14}+\r_{42})/2$.

At last, if the pair $\{p_1,\,p_3\}$ is a mustaches of $G$, then the weight of the corresponding minimal parametric filling is more than or equal $\r_{13}+\r_{24}$, but, due to Condition~$(*)$, this value is more than or equal to the weights of minimal parametric fillings with two previous topologies, therefore, minimal filling is among those two minimal parametric fillings. Comparing the weights of these two fillings, we see that the minimality of the filling with mustaches  $\{p_1,\,p_2\}$ is equivalent to condition $\r_{12}+\r_{43}\le\r_{23}+\r_{14}$. Taking into account Condition~$(*)$, the latter inequality is equivalent to the next relation
 $$
\r_{12}+\r_{34}=\min\{\r_{12}+\r_{34},\,\r_{13}+\r_{24},\,\r_{14}+\r_{23}\},
 $$
and so, we obtain the required formula. The proof is complete.
 \end {proof}

We apply the result obtained to the vertex set of a planar convex quadrangle.

 \begin {cor}\label{cor:conv4gon}
Let $M$ be the vertex set of a convex quadrangle $p_1p_2p_3p_4\ss\R^2$ and $\r(p_i,p_j)=\|p_i-p_j\|$. The weight of a minimal filling of the space $(M,\r)$ is equal to
 $$
\frac12\min\bigl(\r_{12}+\r_{34},\,\r_{14}+\r_{23}\bigr)+
\frac{\r_{13}+\r_{24}}2.
 $$
The topology of minimal filling is a binary tree with mustaches corresponding to opposite sides of the less total length. If these sides are  $p_1p_2$ and $p_3p_4$, then the weights of edges are described by System~$(**)$.
 \end {cor}

Consider several particular examples. Let $M$ be the vertex set of the diamond consisting of two regular triangles with side $1$, and let $p_2p_4$ be its smaller diagonal. Then $\r_{12}=\r_{23}=\r_{34}=\r_{14}=\r_{24}=1$, and $\r_{13}=\sqrt3$. In this case, in accordance with Assertion~\ref{cor:conv4gon}, there are two symmetric minimal fillings with the weight $(3+\sqrt3)/2$. The weights of the filling with mustaches, say, $\{p_1,\,p_2\}$, have the form
 \begin {gather*}
\l_1=\frac{\sqrt3+1}2-\l_3,\quad
\l_2=\l_3-\frac{\sqrt3-1}2,\quad
\frac12\le\l_3\le\frac{\sqrt3}2,\quad
\l_4=1-\l_3, \\
\l_5=\frac{\sqrt3-1}2.
 \end {gather*}

Now let the quadrangle $p_1p_2p_3p_4$ be a rectangle with the sides $a$ and $b$, and let $a\le b$. The length of the diagonal is denoted by $d=\sqrt{a^2+b^2}$. Then the weight of minimal filling is equal to $a+d$, the topology of minimal filling is the binary tree whose mustaches are the pairs of the vertices of the smallest side, and the weights of the edges are as follows:
 \begin {gather*}
\l_1=a-\l_3,\quad
\l_2=\l_3,\quad
\frac{a+b-d}2\le\l_3\le\frac{a-b+d}2,\quad
\l_4=a-\l_3, \\
\l_5=d-a.
 \end {gather*}

 \subsection {Convex Polygons}
O.~Kolcheva and E.~Zaval'nuk, students of mechanical and mathematical faculty of Moscow state university, study minimal fillings of metricspaces formed by the vertex sets of planar convex polygons. As a result, maximal half-perimeters for convex $5$-gon and $6$-gon are found (O,~Kolcheva), and also for regular polygon (E.~Zaval'nuk). These results give an opportunity to obtain low estimates for the weights of minimal fillings in these cases.

 \section {Ratios} \label{sec:ratios}
 \markright {\thesection.~Ratios.}
The Steiner ratio is an important characteristic in Steiner minimal networks theory. Let $M$ be a finite subset of a metric space $\cX$, consisting of more than one point. Recall that the {\em Steiner ratio of $M$} is the ratio of the lengths of Steiner minimal tree and minimal spanning tree constructed on $M$, i.e\. the value
$\sr(M)=\smt(M)/\mst(M)$. The infimum  of the numbers $\sr(M)$ over all such subsets $M$ of $\cX$ is called by the {\em Steiner ratio of the space $\cX$} and is denoted by  $\sr(\cX)$, see~\cite{GilPol}.

Notice that the exact values of the Steiner ratio are known for a very restricted class of spaces (see a review in~\cite{ITLup} or in~\cite{ITRFFIBook}). Below, we define other two ratios based on minimal fillings, which could be more available for calculating, and which could be useful to calculate the Steiner ratio, as we hope.

In this section we give only some preliminaries of the theory of these ratios. In next publications we hope to develop this theory.

 \subsection {Steiner--Gromov Ratio}
For convenience, the sets consisting from more than a single point are referred as {\em non-trivial}. Let $\cX=(X,\r)$ be an arbitrary metric space, and let $M\ss X$ be some finite subset. Recall that by $\mst(M,\r)$ we denote the length of minimal spanning tree of the space $(M,\r)$. Further, for non-trivial $M$, we define the value
 $$
\sgr(M)=\mf(M,\r)/\mst(M,\r)
 $$
and call it by the {\em Steiner--Gromov ratio\/} of the subset $M$; The value $\inf\sgr(M)$, where the infimum is taken over all non-trivial finite subsets of $\cX$, consisting of at most $n>1$ vertices is denoted by $\sgr_n(\cX)$ and is called by the {\em degree $n$ Steiner--Gromov ratio of the space $\cX$}; at last, the value $\inf\sgr_n(\cX)$, where infimum is taken over all positive integers $n>1$ is called by the {\em Steiner--Gromov ratio of the space $\cX$} and is denoted by $\sgr(\cX)$ or by $\sgr(X)$, if it is clear what particular metric  on $X$ is considered. Notice that $\sgr_n(\cX)$ is non-increasing function on $n$. besides, it is easy to see that $\sgr_2(\cX)=1$ for any non-trivial metric space $\cX$.

 \begin {examp}
Let $\cX=(X,\r)$ be a metric space containing a regular triangle. Calculate $\sgr_3(\cX)$. Consider an arbitrary triangle $M=\{p_1,\,p_2,\,p_3\}\ss X$ and put $d_i=\r(p_j,p_k)$, where $\{i,j,k\}=\{1,2,3\}$. The results of Subsection~\ref{subsec:triangle} imply that $\mf(M,\r)=(d_1+d_2+d_3)/2$. Further, without loss of generality, assume that $d_1\le d_2\le d_3$, then $\mst(M,\r)=d_1+d_2\le 2d_3$, and, hence,
 $$
\sgr_3(\cX)=\inf_M\dfrac{d_1+d_2+d_3}{2(d_1+d_2)}\ge
\dfrac{d_1+d_2+(d_1+d_2)/2}{2(d_1+d_2)}=\frac34.
 $$
On the other hand, for a regular triangle $M\ss X$, we have $\sgr(M)=3/4$, thus $\sgr_3(\cX)=3/4$.
 \end {examp}

 \begin {ass}\label{ass:steiner_ratio}
The Steiner--Gromov ratio of an arbitrary metric space is not less than $1/2$. There exist metric spaces, whose Steiner--Gromov ratio equals to $1/2$.
 \end {ass}

 \begin {proof}
At first, we show that $\sgr(\cX)\ge1/2$ for an arbitrary metric space $\cX=(X,\r)$. Let $M$ be a finite subset of $\cX$. Due to Corollary~\ref{cor:fill_ge_halfp}, we have:
 $$
\mf(M,\r)\ge p(M,\r)>\frac12\mst(M,\r),
 $$
and, hence, $\sgr(\cX)=\inf_M\mf(M,\r)/\mst(M,\r)\ge1/2$, where infimum is taken over all finite non-trivial  $M\ss X$.

Now, let us construct a metric space, whose Steiner--Gromov ratio is equal to $1/2$. Let us take a countable set $X$ with the distance function $\r$ which is equal to  $1$ on any pair of distinct points.

Let $M_n$ be an arbitrary subset of $\cX$, consisting of  $n$ points. The space $(M_n,\r)$ is a regular simplex, therefore, the results of Subsection~\ref{subsec:simplex} imply that $\mf(M_n,\r)=n/2$. On the other hand, $\mst(M_n,\r)=n-1$, therefore,
 $$
\sgr(M_n)=\frac12\dfrac{n}{n-1}.
 $$
As $n\to\infty$, the value $\sgr(M_n)$ tends to $1/2$, so $\sgr(\cX)\le1/2$. Taking into account the above inequality, we conclude that $\sgr(\cX)=1/2$. The proof is complete.
 \end {proof}

 \subsection {Steiner Subratio}
Let $\cX=(X,\r)$ be an arbitrary metric space, and let $M\ss\cX$ be some its finite subset. Recall that by $\smt(M,\r)$ we denote the length of Steiner minimal tree joining  $M$. Further, for non-trivial subsets $M$, we define the value
 $$
\ssr(M)=\mf(M,\r)/\smt(M,\r)
 $$
and call it by the {\em Steiner subratio\/} of the set $M$; the value $\inf\ssr(M)$, where infimum is taken over all non-trivial finite subsets of $\cX$, consisting of at most $n>1$ points, is denoted by $\ssr_n(\cX)$ and is called by the {\em degree $n$ Steiner subratio of the space $\cX$}; at last, the value $\inf\ssr_n(\cX)$, where infimum is taken over all positive integers  $n>1$, is called by the {\em Steiner subratio of the space $\cX$} and is denoted by $\ssr(\cX)$ or by $\ssr(X)$, if it is clear what particular metric  on $X$ is considered. Notice that $\ssr_n(\cX)$ is a non-increasing function on $n$.  Besides, it is easy to see that $\ssr_2(\cX)=1$ for any non-trivial metric space $\cX$.

 \begin {examp}
Let us show that $\ssr_3(\R^n,\r)=\ssr_3(\R^n)$, $n\ge2$, where $\r$ is the standard Euclidean metric is equal to $\sqrt3/2$. Let $M\ss\R^n$ be a triangle with the sides $a\le b\le c$, then, in accordance with the results of Subsection~\ref{subsec:triangle}, $\mf(M,\r)=(a+b+c)/2$. For the triangles having an angle which is more than or equal to $120^\c$, we have $\smt(M)=a+b$, and, hence,
 \begin {multline*}
\ssr_3(M)=\frac{a+b+c}{2(a+b)}\ge
\frac{a+b+\sqrt{a^2+b^2+ab}}{2(a+b)}=
\frac12+\frac{\sqrt{(a+b)^2-ab}}{2(a+b)}=\\ =
\frac12+\frac12\sqrt{1-\frac1{(\sqrt{a/b}+\sqrt{b/a})^2}}\ge
\frac12+\frac12\sqrt{1-\frac14}=
\frac12+\frac{\sqrt3}4>\frac{\sqrt{3}}2.
 \end {multline*}

Now, let all angles of the $M$ be less than $120^\c$. Then the shortest tree possesses an additional vertex $s$ adjacent with all the points from $M$, and the edges meet at $s$ by the angles of $120^\c$. Let $M=\{p_1,p_2,p_3\}$, and let $x$, $y$, and $z$ be the lengths of the segments $sm_1$, $sm_2$, and $sm_3$, respectively. Put $f(u,v)=\sqrt{u^2+v^2+uv}$, then the perimeter of the triangle $M$ is equal to $P(x,y,z)=f(x,y)+f(x,z)+f(y,z)$. Since the value $\ssr_3(\R^n)$ is preserved under similarity transformations, we have
 $$
\ssr_3(\R^n)=\inf\frac{P(x,y,z)/2}{x+y+z},
 $$
where $\inf$ is taken over all positive $x$, $y$, and $z$, which $x+y+z$ is constant, say equals $1$. Thus, to calculate $\ssr_3(\R^n)$ it suffices to find the infimum of the function $P(x,y,z)$ under the above restriction and to compare the result with the above estimates.

Extending the function $P(x,y,z)$ to the continuous function defined for zero values of variables, we obtain a continuous function on a closed triangle. Such function $P(x,y,z)$ achieves its minimal value. Show that the minimum is achieved at the point  $(1/3,1/3,1/3)$. i.e\. at the regular triangle.

Indeed, consider some point $(x,y,z)$ from the domain of the function $P(x,y,z)$ and assume that not all the coordinates of this point are equal to each other. Without loss of generality, assume that $x\ne y$. Consider the new point with the coordinates $\bigl((x+y)/2,(x+y)/2,z\bigr)$, the sum of whose coordinates is equal to $1$. Show that the value of $P$ at this point is less than the initial one (that proves that $P$ achieves minimum at the regular triangle $M$). We have:
 $$
P\bigl((x+y)/2,(x+y)/2,z\bigr)=
\frac{\sqrt3}2(x+y)+\sqrt{(x+y)^2+4z^2+2xz+2yz}.
 $$
Show that
 \begin {gather*}
\frac{\sqrt3}2(x+y)<f(x,y)\ \ \text{и}\\
\sqrt{(x+y)^2+4z^2+2xz+2yz}\le f(x,z)+f(y,z).
 \end {gather*}
The first inequality is equivalent to
 $$
0<f^2(x,y)-\left(\frac{\sqrt3}2(x+y)\right)^2=
x^2+y^2+xy-\frac34(x+y)^2=\frac14(x-y)^2,
 $$
which is valid, since $x\ne y$.

Now, let us prove the second inequality. It is equivalent to
 \begin {multline*}
0\le(f(x,z)+f(y,z))^2-\bigl((x+y)^2+4z^2+2xz+2yz\bigr)=\\ =
2f(x,z)f(y,z)-2xy-xz-yz-2z^2,
 \end {multline*}
which, by turn, is equivalent to
 $$
0\le4f^2(x,z)f^2(y,z)-(2xy+xz+yz+2z^2)^2=3(x-y)^2z^2,
 $$
that is valid, also.

Thus, we show that $\ssr_3(\R^n)\ge\sqrt3/2$, and the equality is achieved at the regular triangle. Comparing this result with the infimum taken over the triangles having an angle of more than or equal to  $120^\c$ we conclude that the infimum over all triangles is achieved at the regular triangle. Thus, $\ssr_3(\R^n)=\sqrt3/2$.
 \end {examp}

The next result is obtained by E.~Filonenko, a student of mechanical and mathematical department of Moscow state university, see~\cite{Filonenko}.

 \begin {ass}[E.~Filonenko]
$\ssr_4(\R^2)=\sqrt3/2$.
 \end {ass}

 \begin {conj}\label{conj:subrat}
The Steiner subratio of the Euclidean plane is achieved at the regular triangle and, hence, is equal to $\sqrt3/2$.
 \end {conj}

 \markright{References.}

 \end {document}